\documentclass[final,1p,times,authoryear]{elsarticle}
\pdfoutput=1

\usepackage{amssymb}
\usepackage{amsthm}

\theoremstyle{definition}
\newtheorem{definition}{Definition}[section]

\usepackage[fleqn]{amsmath}
\DeclareMathOperator*{\argmin}{arg\,min}

\usepackage{url}

\usepackage{mathtools}
\DeclarePairedDelimiter{\ceil}{\lceil}{\rceil}

\usepackage{multirow}

\begin{document}

\begin{frontmatter}



\title{The Benefits of Autonomous Vehicles for Community-Based Trip Sharing}


\author[1]{Mohd. Hafiz Hasan}
\ead{hasanm@umich.edu}

\author[2]{Pascal Van Hentenryck\corref{cor1}}
\ead{pascal.vanhentenryck@isye.gatech.edu}

\cortext[cor1]{Corresponding author}

\address[1]{University of Michigan, Ann Arbor, Michigan 48105, USA}
\address[2]{Georgia Institute of Technology, Atlanta, Georgia 30332, USA}

\begin{abstract}

  This work reconsiders the concept of community-based trip sharing
  proposed by \citet{hasan2018} that leverages the structure of
  commuting patterns and urban communities to optimize trip
  sharing. It aims at quantifying the benefits of autonomous vehicles
  for community-based trip sharing, compared to a car-pooling platform
  where vehicles are driven by their owners. In the considered
  problem, each rider specifies a desired arrival time for her inbound trip
  (commuting to work) and a departure time for her outbound trip
  (commuting back home). In addition, her commute time cannot deviate
  too much from the duration of a direct trip.  Prior work motivated
  by reducing parking pressure and congestion in the city of Ann
  Arbor, Michigan, showed that a car-pooling platform for
  community-based trip sharing could reduce the number of vehicles by
  close to 60\%.

  This paper studies the potential benefits of autonomous vehicles in
  further reducing the number of vehicles needed to serve all these
  commuting trips. It proposes a column-generation procedure that
  generates and assembles mini routes to serve inbound and outbound
  trips, using a lexicographic objective that first minimizes the
  required vehicle count and then the total travel distance.  The
  optimization algorithm is evaluated on a large-scale, real-world
  dataset of commute trips from the city of Ann Arbor, Michigan. The
  results of the optimization show that it can leverage
  autonomous vehicles to reduce the daily vehicle usage by 92\%,
  improving upon the results of the original Commute Trip Sharing
  Problem by 34\%, while also reducing daily vehicle miles traveled by
  approximately 30\%. These results demonstrate the significant
  potential of autonomous vehicles for the shared commuting of a
  community to a common work destination.
\end{abstract}



\end{frontmatter}


\section{Introduction} 
\label{sec:intro}

The notion of community-based trip sharing---leveraging the structure
of commuting patterns and urban communities when optimizing trip
sharing---was first explored by \citet{hasan2018} to reduce parking
pressure and congestion on university and corporate campuses. The study,
which was originally motivated by the desire to relieve parking
pressure at the University of Michigan, Ann Arbor, investigated the 
effects of different driver and commuter matching arrangements on trip 
shareability for a car-pooling or car-sharing platform. Trip 
shareability was loosely defined as the ability to aggregate as many trips 
as possible to reduce the number of vehicles required to serve them. 
The evaluation of several different optimization models revealed that commuter
matching flexibility, i.e., a willingness to be matched with
different drivers and passengers daily, is key for an effective
trip-sharing platform.

This early work was extended by \citet{hasan2020}, where the Commute
Trip Sharing Problem (CTSP) was formalized. The CTSP seeks a routing
plan for the set of commuting trips that minimizes a lexicographic
objective. The primary objective is to minimize the number of vehicles
to cover all trips, while the secondary objective is to minimize the
total travel distance. Every commuting trip consists of a pair of trip
requests, one to the workplace (inbound) and another back home
(outbound), each with specific pickup and drop-off locations as well
as time windows specifying allowable service times at each
location. The routes of the CTSP must serve each request exactly once
and ensure that a vehicle capacity and trip specific ride-duration
limits are not exceeded. The ride-duration constraint guarantees a
level of quality of service for the riders. Finally, as the routing
plan selects drivers from the set of commuters, it must ensure that
the set of drivers selected for the inbound trips is
identical to that for the outbound trips.  The CTSP is thus a Vehicle
Routing Problem (VRP) with time-window, capacity, pairing, precedence,
ride-duration, and driver constraints. \citet{hasan2020} proposed two
exact algorithms to solve the CTSP and applied them on a real-world
dataset from the city of Ann Arbor, Michigan. The case study shows
that community-based car pooling can decrease daily vehicle usage by
up to 57\%. These results highlighted the significant potential in
vehicle reduction of community-based trip sharing. However, the
vehicles in the CTSP routing plans are still mostly idle, as they
perform a single inbound and outbound trip a day. Moreover, the
constraint that a trip starts at the driver origin and ends at the
driver destination limits the potential for ridesharing. These limitations
make the adoption of autonomous vehicles (AVs) for the CTSP particularly
appealing, as the absence of drivers would directly address these key shortcomings 
and could potentially lead to further reductions in fleet size, higher 
vehicle utilization, and increased ridesharing.

The goal of this paper is to examine the potential benefits of using
AVs for performing the same trips, quantifying the
reduction in fleet size and the miles travelled. It studies the
Commute Trip Sharing Problem for Autonomous Vehicles (CTSPAV) which is
similar to CTSP but uses a fleet of AVs that depart
from and return to a designated depot to serve all commute trips. The
CTSPAV is very similar to, and is a specialization of, the Dial-A-Ride
Problem (DARP) \citep{cordeau2006}.  The main difference is that trip
requests in the CTSPAV come in pairs, one to the workplace (typically
in the morning) and another to return back home (typically in the
evening). This feature makes it possible to adopt solution techniques
that are computationally attractive. Our proposed solution approach 
involves chaining mini routes, i.e., short routes with distinct pickup, transit,
and drop-off phases, to form longer AV routes. It is especially suited
for problem scenarios involving commuters traveling to a common/centralized 
location, e.g., the commute trips of employees of a university or corporate campus, 
or those involving commuters living in a common/centralized location,
e.g., the commute trips originating from an apartment complex or a residential neighborhood.
The case study considered in this work, which consists of commute trips made
to 15 university-owned parking structures located within close vicinity
to each other in downtown Ann Arbor, Michigan, fits the profile of the former scenario,
making it a prime candidate for evaluating the efficacy of the
proposed approach.   

The main contributions of this paper can be summarized as follows:
\begin{enumerate}[1.]
\item The paper formalizes the CTSPAV that seeks an optimal set of routes for
  a fleet of AVs for serving a set of commute trips subject to the
  passenger-related constraints of the original CTSP.

\item The paper proposes a column-generation procedure to find a
  high-quality solution to the problem. The procedure uses a pricing
  problem to generate mini routes that are then assembled in a master
  problem. Each mini route serves inbound (resp. outbound) trips for a
  number of riders, satisfying time-window and ride-duration
  constraints on the trips as well as the capacity constraints on the
  vehicles. This approach is a departure from classical column-generation 
  procedures commonly adopted for VRPs and the DARP, whereby the pricing problem
  searches for complete routes departing from and returning to
  a depot, while the master problem solves a set-covering/partitioning 
  problem that ensures every customer is served.

\item The paper shows that the  proposed algorithm outperforms a
  state-of-the-art DARP algorithm based on the
  classical column-generation approach for the CTSPAV.

\item The paper applies the proposed algorithms on a large-scale,
  real-world dataset of commute trips for the city of Ann Arbor,
  Michigan. The experimental results show that the algorithm is capable of 
  reducing the daily vehicle usage by 92\%, improving upon the results of the
  original CTSP by 34\%, while also reducing
  daily vehicle miles traveled by approximately 30\%.
\end{enumerate}

\noindent
Overall, the results demonstrate the significant potential of AVs for
serving the commuting needs of a community whose members work at a
common location. The paper also includes a coarse cost analysis that
highlights that fleet sizing is the correct metric to optimize when the goal consists of maximizing
the profitability of the service. 

The rest of this paper is organized as follows. Section \ref{sec:related} briefly 
outlines some related work while Section
\ref{sec:notation} presents the terminology and notations used throughout
the paper. Section \ref{sec:CTSPAV} specifies the CTSPAV and
presents a mixed-integer programming (MIP) model that 
formalizes the CTSPAV. Section
\ref{sec:CG} describes the column-generation procedure. Section
\ref{sec:darp} sketches a DARP-based procedure that is used for comparison
purposes. Section \ref{sec:results} presents detailed computational
results and examines the performance of the proposed approach for
servicing the commuting needs of the community under study.  Finally,
Section \ref{sec:conclusion} provides some concluding remarks.

\section{Related Work} 
\label{sec:related}

Both the CTSP and the CTSPAV are generalizations of the Vehicle Routing Problem with Time Windows (VRPTW). The VRPTW seeks a set of minimum cost routes, each departing from and returning to a designated depot, that serve a set of customers, each with a capacity demand and a time window within which service may commence. The problem must ensure that each customer is served exactly once within their time windows while not allowing the vehicle capacities to be exceeded. It is well known to be NP-hard, as \cite{savelsbergh1985} showed that finding a feasible solution to the problem for a fixed vehicle count is strongly NP-complete. It has been studied extensively in the literature, and numerous methods have been suggested to tackle its complexity, from metaheuristics like \cite{taillard1997} and \cite{braysy2005}, to exact solution methods using Lagrangian relaxation \citep{kohl1997,kallehauge2006}, column generation \citep{desrosiers1984,desrochers1992}, or polyhedral approaches \citep{kohl1999,bard2002,kallehauge2007}. See \cite{cordeau2002} for an extensive review of the topic.

The VRPTW was generalized to the Pickup and Delivery Problem with Time Windows (PDPTW) by \cite{dumas1991}, whereby service locations come in pairs, a pickup and a delivery location for each customer, that must be serviced in order by the same route. They proposed a dynamic-programming, label-setting algorithm to search for routes that satisfy the new pairing and precedence constraints along with the existing time-window and vehicle-capacity constraints. The algorithm is incorporated within a column-generation procedure to solve the problem. A similar approach was also adopted by \cite{ropke2009} in their branch-and-cut-and-price algorithm for the PDPTW, while \cite{ruland1997} used a polyhedral approach to solve the version of the problem without time-window constraints. The DARP generalizes the PDPTW by introducing ride-duration constraints which limit the maximum duration between each pickup and delivery location pair. It models the maximum time spent on the vehicle by every customer and is critical for guaranteeing a quality-of-service level for services that transport passengers instead of merchandise, like door-to-door transportation services for the disabled and the elderly or those for ridesharing. The problem has also been extensively reviewed by \cite{cordeau2003a} and \cite{cordeau2007}, and it has been tackled with methods ranging from approximate methods like heuristics \citep{bodin1986, jaw1986} and metaheuristics \citep{cordeau2003, ritzinger2016}, to exact ones utilizing cutting plane methods \citep{cordeau2006} and column generation \citep{gschwind2015}. 

Of the many solution approaches proposed for the various generalizations of the VRPTW, column generation is perhaps the most popular due to its elegance in only considering routes that can improve the objective function, and its proven effectiveness in producing strong lower bounds to the problem objective when used in conjunction with the Dantzig--Wolfe decomposition \citep{dantzig1960}. It entails a pricing subproblem which searches for routes satisfying problem-specific feasibility constraints, a problem that is typically an Elementary Shortest Path Problem with Resource Constraints (ESPPRC), whereby resource contraints are used to model the feasibility constraints, and the elementarity requirement ensures that each customer is serviced exactly once. Unfortunately, the ESPPRC has been proven to be NP-hard in the strong sense by \cite{dror1994}, and while exact solution methods have been proposed for the problem, e.g., \cite{feillet2004,chabrier2006,boland2006,drexl2013}, the elementarity requirement is more commonly relaxed to produce a Shortest Path Problem with Resource Constraints (SPPRC) which admits a pseudo-polynomial solution approach. A variety of strategies are then adopted to handle non-elementary paths, e.g., \cite{desrosiers1984}, \cite{dumas1991}, \cite{ropke2009}, and \cite{gschwind2015} eliminate them by either preventing their selection in an integer solution or by using infeasible path elimination constraints in the master problem, while \cite{desrochers1992} and \cite{irnich2006} take a middle-ground approach by eliminating 2- and $k$-cycles from discovered paths respectively. Regardless of whether an SPPRC or an ESPPRC is used in the pricing subproblem, they are typically solved via dynamic programming, the most popular being the generalized label-setting algorithm for multiple resource constraints by \cite{desrochers1988b}. Other suggested dynamic-programming approaches include the label-correcting algorithm by \cite{desrosiers1983} which is based on the Ford-Bellman-Moore algorithm and the label-setting algorithm by \cite{desrochers1988} which generalizes Dijkstra's algorithm. Methods utilizing Lagrangian relaxation \citep{beasley1989,borndrfer2001}, constraint programming \citep{rousseau2004}, and cutting planes \citep{drexl2013} have also been proposed, and \cite{irnich2005} provides an in-depth review  of the SPPRC.

More recently, an increased awareness for sustainability of passenger transportation systems combined with the availability of large-scale, real-world trip datasets has shifted the focus towards optimization of car-pooling and ride-sharing services to reduce traffic congestion and pollution. \cite{baldacci2004} studied the Car-Pooling Problem (CPP) which seeks to minimize the number of private cars used for commuting to a common workplace. They considered a variant of the problem which optimizes car pooling for the trips to the workplace independently from those for the return trips and assumes that the set of drivers and passengers are known beforehand, making it a specialization of the DARP. The effectiveness of their proposed Lagrangian column-generation method was demonstrated on instances derived from real-world data provided by a research institution in Italy. \cite{agatz2011} contrasts the CPP with the dynamic ride-sharing problem, whereby the latter matches drivers and riders for single, non-recurring trips in real time. They proposed an optimization method which casts the problem as a graph matching problem and solves it at regular intervals within a rolling-horizon framework. They also presented a case study which applies the approach on real-world travel demand data from metro Atlanta. \cite{santi2014} introduced the notion of shareability graphs as a tool to quantify the potential benefits of ridesharing, and applied it on trip data from the New York City (NYC) Taxi and Limousine Commission trip record which stores information of more than one billion taxi rides in NYC recorded since January 2009. \cite{alonso-mora2017} then built on this idea to propose an anytime optimal algorithm for the on-demand ride-sharing problem, and the efficacy of their method was also demonstrated through its application on the trips from the NYC taxi dataset. \cite{agatz2012} discusses the different planning considerations for and the issues arising from dynamic ridesharing by classifying the different variations of the problem and reviewing the optimization approaches proposed for them. \cite{mourad2019} takes a broader view of shared mobility in their survey, whereby applications which combine transportation of people and freight in both pre-arranged and real-time settings are reviewed together with their corresponding optimization approaches.

The advent of self-driving technology combined with the race to achieve full driving automation has also triggered a growing interest in Shared Autonomous Vehicle (SAV) systems. Advances in Shared Mobility Services (SMS) and AV technology are widely considered to be mutually beneficial, as the widespread adoption of AVs in SMS could help make AVs financially viable \citep{gurumurthy2018,stocker2019} and accelerate the proliferation of SMS \citep{thomas2018} at the same time. The potential impact of SAV services, ranging from their effect on the economy and the environment to the changes in policy necessary for their governance, have also been widely discussed and reviewed in works like \cite{milakis2017}, \cite{soteropoulos2019}, and \cite{narayanan2020}. \cite{narayanan2020} also proposed classifying SAV services as either on-demand or reservation-based systems, according to the time frame within which the trip requests are made. The former allows requests to be made in real time, making it better suited for serving dynamic trips, whereas the latter requires requests to be made in advance, making it better for recurring trips. Each has its own set of advantages: While on-demand systems can address dynamically changing trip demand, reservation-based systems can further reduce the fleet size and increase the efficiency of routes (by reducing empty cruising time and increasing the number of customers served per vehicle), as demonstrated by \cite{wang2014}, as they know the requests ahead of time and can optimize trips over a longer time horizon. Several optimization approaches have been proposed for both systems. For on-demand systems, \cite{farhan2018} proposed a three-step approach for optimizing a fleet of SAVs that serves on-demand trips. It first discretizes the time horizon into 5-minute intervals, clusters trip requests from each interval by assigning the riders to their nearest vehicle, and finds the optimal vehicle routes by modeling the requests from each cluster as a VRPTW and solving the problem using a tabu-search metaheuristic. On the other hand, \cite{pinto2020} considered integrating SAVs with an existing public transit system to better serve lower density areas using a bi-level modeling framework to jointly optimize the transit network schedule together with the sizing of the AV fleet. They proposed an iterative heuristic which solves a transit network frequency setting problem using a non-linear solver in the upper level and solves a dynamic combined mode choice-traveler assignment problem using an agent-based simulation in the lower level. For reservation-based systems (which are similar to the system proposed in this study), \cite{ma2017} proposed an approach to optimize a fleet of SAVs for trips requests that are known ahead of time. However, their approach only allows vehicle sharing whereby each trip is served without being interrupted by other trip requests. This restriction admits an LP model for the problem which can then be solved efficiently. The modeling technique, however, is not applicable to ridesharing problems like to one considered in this study. \cite{bongiovanni2019} considers a variant of the DARP that uses electric autonomous vehicles, called the e-ADARP. It extends the classical DARP by incorporating additional considerations, like battery management and intermediate stops for vehicle recharging, that only apply to the operation of electric AVs. They proposed two- and three-index formulations for the problem which are solved using a traditional branch-and-cut approach which incorporates new, problem specific valid inequalities. They demonstrated the approach's ability to produce optimal solutions for instances with up to 40 trip requests. Our study, however, considers instances that are five times larger and would therefore require a more robust approach. Numerous other works have touted the potential benefits of these SAV systems, from reducing traffic \citep{martinez2017,alazzawi2018,salazar2018}, to increasing road capacity \citep{friedrich2015,tientrakool2011,talebpour2016,menaoreja2018,olia2018}, to reducing parking demand \citep{zhang2015,dia2017,zhang2017}. However, there also appears to be a consensus that the benefits require AV adoption reaching a critical mass before they can be truly realized.

\cite{hasan2018} introduced community-based trip sharing, in which they investigated the efficacy of several optimization models which utilize different sets of driver and passenger matching conditions to optimize car pooling and car sharing for commuting and considered a community-based partitioning approach which clusters commuters based on their residential neighborhoods. They concluded that commuter matching flexibility, i.e., their willingness to adopt different roles and to be matched with different drivers and passengers daily, is critical for an effective trip-sharing platform. This work was extended by \cite{hasan2020}, whereby the best performing car-pooling model from their previous work, which selects the optimal set of drivers from the set of commuters and optimizes their inbound and outbound routes on a daily basis to reduce vehicle count and their total travel distance, is formalized as the CTSP. They then proposed two exact solution approaches for the problem: A method which first searches for all feasible routes and then optimizes their selection using an integer program, and a branch-and-price algorithm which searches for feasible routes using column generation. Evaluations on a large-scale, real-world commute-trip dataset from the city of Ann Arbor, Michigan, revealed their capability to reduce vehicle count by up to 57\%. \cite{hasan2020b} then considered special variants of the CTSP in which uncertainties are associated with the schedules of return trips, and proposed a stochastic optimization approach which uses scenario sampling to effectively solve versions of the problem whereby commuters confirm their return times by a fixed deadline or in real time. This paper considers another natural extension to the CTSP. It specifically aims at investigating the potential of autonomous vehicles to address a key shortcoming of the CTSP: the short routes induced by its driver constraints that limit the effective use of the vehicles. Access to the Ann Arbor commute trip dataset and the results of the original CTSP on the trips from the dataset allows this work to be uniquely poised to compare and contrast the efficacy of different optimization approaches for conventional and autonomous vehicles, and to provide new insights into the benefits afforded by autonomous vehicles for large-scale sharing of commute trips.

\section{Preliminaries} 
\label{sec:notation}

This section defines the main concepts used in this paper: trips, mini
routes, and AV routes. It also describes the constraints that mini
routes and AV routes must satisfy. This work assumes the utilization of a
homogeneous fleet of vehicles with capacity $K$ to serve all rides,
and that the triangle inequality is satisfied for all travel times and distances.

\paragraph{Trips} A trip $t=\{o,dt,d,at\}$ is a tuple that consists of an origin $o$, 
a departure time $dt$, a destination $d$, and an arrival time $at$ of a trip request. 
Every day, a commuter $c$ makes two trips: a trip $t_c^+$ to the workplace and a
return trip $t_c^-$ back home. These trips are called inbound and
outbound trips respectively. 

\paragraph{Mini Routes} A mini route $r$ is a sequence of locations that visits each origin
and destination from a set of inbound or outbound trips exactly
once. Let $\mathcal{C}_r$ denote the set of riders served in $r$. A
mini route $r$ must respect the vehicle capacity, i.e.,
$|\mathcal{C}_r| \leq K$, and consists of three phases: a pickup phase
where the passengers are picked up, a transit phase where the vehicle
travels to the destination, and a drop-off phase where all the
passengers are dropped off. During the pickup (resp. drop-off) phase,
the vehicle visits only origins (resp. destinations), whereas it
travels from an origin to a destination in the transit phase. For
instance, a possible mini route for a car with $K=4$ serving trips
$t_1=\{o_1,dt_1,d_1,at_1\}$, $t_2=\{o_2,dt_2,d_2,at_2\}$, and
$t_3=\{o_3,dt_3,d_3,at_3\}$ is $r=o_2\rightarrow o_1\rightarrow
o_3\rightarrow d_1\rightarrow d_2\rightarrow d_3$, and its pickup,
transit, and drop-off phases are given by $o_2\rightarrow
o_1\rightarrow o_3$, $o_3\rightarrow d_1$, and $d_1\rightarrow
d_2\rightarrow d_3$ respectively. An inbound mini route $r^+$ covers
only inbound trips and an outbound mini route $r^-$ covers only
outbound trips.

\begin{definition}[Valid Mini Route]
	A valid mini route $r$ serving a set $\mathcal{C}_r$ of riders
        visits all of its origins, $\{o_c:c\in\mathcal{C}_r\}$, before
        its destinations, $\{d_c:c\in\mathcal{C}_r\}$, and respects
        the vehicle capacity, i.e., it has $|\mathcal{C}_r|\leq K$.
\end{definition}

\noindent
Let $T_i$ denote the time at which service begins at location $i$,
$s_i$ the service duration at $i$, $pred(i)$ the location visited just
before $i$, $\tau_{(i,j)}$ the estimated travel time for the shortest
path between locations $i$ and $j$, and $\dot{\mathcal{C}}_r$ the
first commuter served on $r$.  Commuters sharing rides are willing to
tolerate some inconvenience in terms of deviations to their desired
departure and arrival times, as well as in terms of their ride
durations compared to their individual, direct trips. Therefore, a
time window $[a_i,b_i]$ is constructed around the desired times and is
associated with each pickup location $i$, where $a_i$ and $b_i$ denote
the earliest and latest times at which service may begin at $i$
respectively. Conversely, only an upper bound $b_j$ is associated with
each drop-off location $j$ as the arrival time at $j$ is implicitly
bounded from below by $a_j = a_i + s_i + \tau_{(i,j)}$, where $i$ is
the corresponding pickup location for $j$. On top of that, a duration
limit $L_c$ is associated with each rider $c$ to denote her maximum
ride duration.

\begin{definition}[Feasible Mini Route]
A feasible mini route $r$ is valid, has pickup and drop-off times
$T_i\in [a_i,b_i]$ for each location $i\in r$, and ensures the ride
duration of each rider $c\in \mathcal{C}_r$ does not exceed $L_c$.
\end{definition}

\noindent
Determining if a valid mini route $r$ is feasible amounts to solving
a feasibility problem defined by the following constraints on $r$.
\begin{flalign}
& a_{o_c}\leq T_{o_c} \leq b_{o_c} \qquad \forall c\in \mathcal{C}_r  \label{eqn:time_window_origin} \\
& T_{d_c} \leq b_{d_c} \qquad \forall c\in \mathcal{C}_r \label{eqn:time_window_destination} \\
& T_{pred(o_c)} + s_{pred(o_c)} + \tau_{(pred(o_c),o_c)} \leq T_{o_c}  \qquad \forall c\in \mathcal{C}_r\setminus\dot{\mathcal{C}}_r \label{eqn:travel_time_origin} \\
& T_{pred(d_c)} + s_{pred(d_c)} + \tau_{(pred(d_c),d_c)} = T_{d_c}  \qquad \forall c\in \mathcal{C}_r \label{eqn:travel_time_destination} \\
& T_{d_c} - (T_{o_c} + s_{o_c}) \leq L_c \qquad \forall c\in \mathcal{C}_r \label{eqn:ride_duration_limit} 
\end{flalign}

\noindent
Constraints \eqref{eqn:time_window_origin} and
\eqref{eqn:time_window_destination} are time-window constraints for
pickup and drop-off locations respectively, while constraints
\eqref{eqn:travel_time_origin} and \eqref{eqn:travel_time_destination}
describe compatibility requirements between pickup/drop-off times and
travel times between consecutive locations along the route. Finally,
constraints \eqref{eqn:ride_duration_limit} specify the ride-duration
limit for each rider. Note that constraints
\eqref{eqn:travel_time_origin} allow waiting at pickup locations.
Moreover, the service starting times on consecutive locations along $r$ 
are strictly increasing, which ensures that the route is
elementary.  Numerous algorithms have been proposed for solving this
feasibility problem efficiently, e.g. \citet{tang2010},
\citet{haughland2010}, \citet{firat2011}, and \citet{gschwind2015}. In
the following, the Boolean function $feasible(r)$ is used to indicate
whether mini route $r$ admits a feasible solution to constraints
\eqref{eqn:time_window_origin}--\eqref{eqn:ride_duration_limit}.

\paragraph{AV Routes} An AV route $\rho=v_s\rightarrow r_1\rightarrow \ldots\rightarrow
r_k\rightarrow v_t$ is a sequence of $k$ distinct mini routes that
starts at a source node $v_s$ and ends at a sink node $v_t$, both
representing a designated depot.

\begin{definition}[Feasible AV Route]
A feasible AV route $\rho$ is one that consists of a sequence of
distinct, feasible mini routes and starts and ends at a designated
depot.
\end{definition}

\noindent
In other words, for $\rho$ to be feasible, each of its mini routes
must be valid and satisfy constraints
\eqref{eqn:time_window_origin}--\eqref{eqn:ride_duration_limit}. Let
$\dot{r}$ denote the first location visited on $r$ and $\ddot{r}$
denote the last. Moreover, each mini route $r_i \; (1 \leq i \leq k)$
must satisfy the following constraints:
\begin{flalign}
& T_{v_s} + \tau_{(v_s, \dot{r}_1)} = T_{\dot{r}_1}  \label{eqn:travel_time_begin} \\
& T_{\ddot{r}_i} + s_{\ddot{r}_i} + \tau_{(\ddot{r}_i, \dot{r}_{i+1})} \leq T_{\dot{r}_{i+1}} \qquad \forall i = 1, \ldots, k-1 \label{eqn:travel_time_routes} \\
& T_{\ddot{r}_k} + s_{\ddot{r}_k} + \tau_{(\ddot{r}_k, v_t)} = T_{v_t} \label{eqn:travel_time_end} 
\end{flalign}

\noindent
Constraints \eqref{eqn:travel_time_begin}--\eqref{eqn:travel_time_end}
describe compatibility requirements between the beginning/ending
service times of consecutive mini routes along $\rho$ and the travel
times between them. The constraints, together with
\eqref{eqn:travel_time_origin} and
\eqref{eqn:travel_time_destination}, enforce strictly increasing
starting times for service on all consecutive locations along $\rho$,
therefore ensuring that $\rho$ is elementary.

\section{The Commute Trip Sharing Problem for Autonomous Vehicles}
\label{sec:CTSPAV}

This section specifies the CTSPAV that seeks a set of AV routes of minimal cost
to serve each inbound and outbound trip of a set of commuters
$\mathcal{C}$ exactly once. Let $n = |\mathcal{C}|$ denote the total
number of commuters, $\mathcal{P}^+ = \{1,\ldots,n\}$ and
$\mathcal{D}^+ = \{n+1,\ldots,2n\}$ denote the sets of all pickup and
drop-off nodes of inbound trips respectively, and $\mathcal{P}^- =
\{2n+1,\ldots,3n\}$ and $\mathcal{D}^- = \{3n+1,\ldots,4n\}$ denote
the sets of all pickup and drop-off nodes of outbound trips
respectively. Let $\mathcal{P} = \mathcal{P}^+\cup\mathcal{P}^-$ and
$\mathcal{D} = \mathcal{D}^+\cup\mathcal{D}^-$.  The nodes have been
defined such that the inbound pickup, inbound drop-off, outbound
pickup, and outbound drop-off locations of commuter $i$ are
represented by nodes $i$, $n+i$, $2n+i$, and $3n+i$ respectively, and
$n+i$ gives the corresponding drop-off node of pickup node
$i\in\mathcal{P}$.

Let $\mathcal{G} = (\mathcal{N}, \mathcal{A})$ denote a directed graph with 
the node set $\mathcal{N} = \mathcal{P}\cup\mathcal{D}\cup\{v_s,v_t\}$ containing 
all pickup and drop-off nodes together with a source and a sink node representing
a designated depot. A ride-duration limit $L_i$ is associated with each node
$i\in\mathcal{P}$. A time window $[a_i,b_i]$ and service duration
$s_i$ are also associated with each node $i\in\mathcal{P}\cup\mathcal{D}$.
There are no time-window constraints for the start and end times of
any AV route, as it is assumed that the AVs may start and end their
routes at any time of the day. In a first approximation, the edge set
$\mathcal{A}=\{(i,j):i,j\in\mathcal{N},i\neq j\}$ consists of all
possible edges. A travel time $\tau_{(i,j)}$, a
distance $\sigma_{(i,j)}$, and a cost $c_{(i,j)}$ are associated with
each edge $(i,j)\in\mathcal{A}$.  The sets of all outgoing and
incoming edges of node $i$ are denoted by $\delta^+(i)$ and
$\delta^-(i)$ respectively.  By definition of AV routes, the following
precedence constraints apply to the set of nodes:
\begin{equation} 
i \prec n + i \prec 2n + i \prec 3n + i \qquad \forall i\in\mathcal{P}^+ \label{eqn:precedence}
\end{equation}
where $i \prec j$ denotes the precedence relation between nodes $i$
and $j$, i.e., the constraint indicating that $i$ must be visited
before $j$ on an AV route.

This paper considers two distinct optimization objectives: (1) a
lexicographic objective that first minimizes the number of vehicles
and then their total travel distance, and (2) a single objective that
only minimizes total travel distance. 

\subsection{A MIP Model for the CTSPAV}
\label{sec:MIP}

This section presents a MIP model for the CTSPAV. The MIP formalizes the CTSPAV and is the foundation of the
column-generation procedure presented in the next section. It is defined in terms of the set $\Omega$ of all 
feasible mini routes and the graph $\mathcal{G}$.

The MIP, also referred to as the master problem
(MP\textsubscript{CTSPAV}), is shown in Figure \ref{fig:CTSPAV}.  It
uses two sets of binary variables: variable $X_r$ indicates whether
mini route $r\in\Omega$ is selected and variable $Y_e$ indicates
whether edge $e\in\mathcal{A}$ is used in the optimal routing plan. It
also uses a continuous variable $T_i$ to represent the start of
service time at node $i\in\mathcal{P}\cup\mathcal{D}$. The model minimizes
the total cost of all selected edges. Constraints
\eqref{eqn:ctspav_route_cover} enforce coverage of each trip by
exactly one mini route, while constraints
\eqref{eqn:ctspav_edge_select} ensure that edges belonging to selected
mini routes are selected. Constraints \eqref{eqn:ctspav_outgoing} and
\eqref{eqn:ctspav_incoming} conserve flow through each pickup and
drop-off node while ensuring each is visited exactly once. Constraints
\eqref{eqn:ctspav_travel_time_1} and \eqref{eqn:ctspav_travel_time_2}
enforce compatibility of service start times and travel times along
selected edges by utilizing large constants for $M_{(i,j)}$ and
$\bar{M}_{(i,j)}$. Constraints \eqref{eqn:ctspav_ride_duration}
describe ride-duration limits for each trip, while constraints
\eqref{eqn:ctspav_time_window_1} are time-window constraints for all
pickup and drop-off nodes.

\begin{figure}[!t]
\begin{flalign}
& \min \sum_{e\in\mathcal{A}} c_e Y_e \label{eqn:ctspav_obj} \\
& \text{subject to} \nonumber \\
& \sum_{r\in\Omega:i\in r} X_r = 1 \qquad \forall i\in\mathcal{P} \label{eqn:ctspav_route_cover} \\
& \sum_{r\in\Omega:e\in r} X_r - Y_e \leq 0 \qquad \forall e\in\mathcal{A} \setminus \{\delta^+(v_s)\cup\delta^-(v_t)\} \label{eqn:ctspav_edge_select} \\
& \sum_{e\in\delta^+(i)} Y_e = 1 \qquad \forall i\in\mathcal{P}\cup\mathcal{D} \label{eqn:ctspav_outgoing} \\
& \sum_{e\in\delta^-(i)} Y_e = 1 \qquad \forall i\in\mathcal{P}\cup\mathcal{D} \label{eqn:ctspav_incoming} \\
& T_i + s_i + \tau_{(i,j)} \leq T_j + M_{(i,j)}(1 - Y_{(i,j)}) \qquad \forall i,j\in\mathcal{P}\cup\mathcal{D} \label{eqn:ctspav_travel_time_1} \\
& T_i + s_i + \tau_{(i,j)} \geq T_j - \bar{M}_{(i,j)}(1 - Y_{(i,j)}) \qquad \forall i\in\mathcal{P}\cup\mathcal{D}, \forall j\in\mathcal{D} \label{eqn:ctspav_travel_time_2} \\
& T_{i+n} - (T_i + s_i) \leq L_i \qquad \forall i\in\mathcal{P} \label{eqn:ctspav_ride_duration} \\
& a_i \leq T_i \leq b_i \qquad \forall i\in\mathcal{P}\cup\mathcal{D} \label{eqn:ctspav_time_window_1} \\
& X_r \in\{0,1\} \qquad \forall r\in\Omega \label{eqn:ctspav_route_var} \\
& Y_e \in\{0,1\} \qquad \forall e\in\mathcal{A} \label{eqn:ctspav_edge_var} 
\end{flalign}
\caption{The CTSPAV Model.}
\label{fig:CTSPAV}
\end{figure}

The lexicographic objective is accomplished using a blended approach
that appropriately weights the sub-objectives: it assigns an
identical, large fixed cost to each AV route and a variable cost that
is proportional to its total distance. Let $\mathcal{R}$ denote the
set of all feasible AV routes. The edge costs are then defined as
follows:
\begin{equation}
c_e=
\begin{cases} \label{eqn:edge_costs}
\sigma_e + 100\cdot\hat{\varsigma}_\text{max}&\qquad\forall e\in\delta^+(v_s)\\
\sigma_e&\qquad\text{otherwise}
\end{cases}
\end{equation}
where $\hat{\varsigma}_\text{max}$ is a constant equal to the length (total distance) of the longest AV route, i.e.:
\begin{equation}
\hat{\varsigma}_\text{max} = \max_{\rho\in\mathcal{R}} \sum_{(i,j)\in\rho} \sigma_{(i,j)}
\end{equation}
The fixed cost, $100\cdot\hat{\varsigma}_\text{max}$, that is significantly larger than the total length of any AV route and that is also identical for every edge $e\in\delta^+(v_s)$ drives the model to first minimize the total flow emanating from the depot, $v_s$. This flow is identical to the total number of AV routes, and thus the number of vehicles used in the solution. The variable cost $\sigma_e$ of every edge $e$  then drives the model to minimize the total distance of all selected edges. The edge costs defined in \eqref{eqn:edge_costs} therefore accomplish the desired lexicographic ordering of the objective, which first minimizes the number of vehicles used in the solution and then their total travel distance.

Conversely, when the objective is to just minimize the total distance, the edge costs are simply
defined by the edge distance, i.e.,
\begin{equation}
c_e = \sigma_e\qquad\forall e\in\mathcal{A}
\end{equation}

The MP\textsubscript{CTSPAV} model can be seen as a scheduling problem
that selects and assembles feasible mini routes to form longer,
feasible AV routes that minimize the total cost. The optimal AV routes
are obtained by constructing paths beginning at $v_s$ and ending at
$v_t$.  The start and end times (at the depot) of these routes are
then obtained by applying equations \eqref{eqn:travel_time_begin} and
\eqref{eqn:travel_time_end} respectively. 

\section{A Column-Generation Procedure for the CTSPAV}
\label{sec:CG}

This section presents a column-generation approach to find
high-quality solutions to the CTSPAV, referred to as the CTSPAV procedure. 
The column-generation approach builds on the MP\textsubscript{CTSPAV} but addresses its main
computational difficulty: the fact that the MP\textsubscript{CTSPAV}
assumes that all mini routes have been pre-computed. The
column-generation approach implements an iterative process that
considers, at each iteration, a subset $\Omega'\subseteq\Omega$ of
feasible mini routes and solves a restricted master problem, denoted
by RMP\textsubscript{CTSPAV}, that is defined as the linear relaxation of
MP\textsubscript{CTSPAV} over $\Omega'$. Using dual information from
RMP\textsubscript{CTSPAV}, the column-generation algorithm then
searches for feasible mini routes with negative reduced costs by
solving a pricing subproblem (PSP\textsubscript{CTSPAV}).  If such
mini routes exist, they are added to the restricted master problem.
Each iteration thus defines a new restricted master problem over a
larger subset of feasible mini routes. The solving of restricted
master problems and pricing subproblems is repeated until the
pricing subproblem cannot find any feasible mini routes with a
negative reduced cost. Upon completion, the optimal objective value of
RMP\textsubscript{CTSPAV} converges to
$z^*$, the optimal objective of the linear relaxation of
MP\textsubscript{CTSPAV}. Whenever the solution of
RMP\textsubscript{CTSPAV} is integral at convergence, it is also
optimal for MP\textsubscript{CTSPAV}. Otherwise, the column-generation
approach solves the final restricted master problem as a MIP
to obtain an integer solution. The objective value of
the MIP provides an upper bound to the optimal solution, while the
objective value $z^*$ of its linear relaxation provides a lower
bound. Together, they are used to compute an optimality gap for the
integer solution.

The approach adopted in this procedure, which has its pricing subproblem search for feasible mini routes which are then chained together in the master problem to form longer, feasible AV routes that depart from and return to the depot, is unlike the classical column-generation approach adopted in most literature on the VRPTW (e.g. \cite{desrosiers1984,desrochers1992}), PDPTW (e.g. \cite{dumas1991,ropke2009}), or DARP (e.g. \cite{gschwind2015})). The classical approach is an application of the Dantzig-Wolfe decomposition on an edge-flow formulation of the problem: it produces a set-partitioning/covering master problem that just selects feasible routes from a set to ensure every customer is served in the solution. Its pricing problem is then solely responsible for searching for \emph{complete} feasible routes that originate from and return to the depot (i.e., routes that would correspond to the AV routes of the CTSPAV). In contrast, our approach shifts part of the burden of constructing the AV routes, which are anticipated to be very long, to the master problem instead of completely relegating the task to the pricing subproblem as is commonly done in the classical approach.

\subsection{The Pricing Subproblem}
\label{sec:PSP}

The PSP\textsubscript{CTSPAV} identifies feasible mini routes with
negative reduced costs. Let $\{\pi_i : i\in\mathcal{P}\}$ and $\{\mu_e
: e\in\mathcal{A}\setminus\{\delta^+(v_s)\cup\delta^-(v_t)\}\}$ denote
dual values associated with constraints \eqref{eqn:ctspav_route_cover}
and \eqref{eqn:ctspav_edge_select} at optimality of
RMP\textsubscript{CTSPAV}.  The reduced cost of mini route $r$ is then
given by:
\begin{equation}\label{mini_route_rc}
\bar{c}_r = -\sum_{i\in r:i\in\mathcal{P}} \pi_i - \sum_{e\in r} \mu_e
\end{equation}

\noindent
The column-generation approach attempts to generate multiple feasible
mini routes during each iteration, one for each pickup node. More
precisely, the PSP\textsubscript{CTSPAV} considers each node
$i\in\mathcal{P}^+\cup\mathcal{P}^-$ as the starting point of a mini
route. For each $i\in\mathcal{P}^+\cup\mathcal{P}^-$, it searches for
a mini route $r_i$ with minimal reduced cost and selects those with
negative reduced costs to augment $\Omega'$. It accomplishes this by
first constructing $2n$ graphs, $\mathcal{G}_i^+$ $(i\in\mathcal{P}^+)$
and $\mathcal{G}_i^-$ $(i\in\mathcal{P}^-)$. It then searches for a least-cost
path from $i$ to a designated sink node that satisfies all mini-route feasibility 
constraints from each graph. The complete details of
this procedure are given in \ref{appendix:pricing}.

\subsection{Practical Implementation Considerations}

This subsection reviews a number of important implementation
techniques for the CTSPAV column-generation procedure.

\paragraph{Filtering of Graph $\mathcal{G}$}

Many edges in $\mathcal{G}$ do not belong to any feasible AV route
and can be removed from $\mathcal{A}$. The following sets of
infeasible edges are obtained by pre-processing time-window, pairing,
precedence, and ride-duration limit constraints on $\mathcal{A}$ using
a combination of rules proposed by \citet{dumas1991} and
\citet{cordeau2006}:
\begin{enumerate}[(a)]
\item Direct trips to and from the depot:
\begin{itemize}
	\item $\{(v_s,v_t), (v_t,v_s)\}$
	\item $\{(i,v_s), (i,v_t), (v_t,i) : i\in\mathcal{P}\}$
	\item $\{(v_s,i), (i,v_s), (v_t,i) : i\in\mathcal{D}\}$
\end{itemize}
\item Pairing and precedence of pickup and drop-off nodes of inbound and outbound trips of each commuter (constraints \eqref{eqn:precedence}): $\{(i,2n+i), (i,3n+i), (n+i,i), (n+i,3n+i), (2n+i,i), (2n+i,n+i), (3n+i,i), (3n+i,n+i), (3n+i,2n+i) : i\in\mathcal{P}^+\}$
\item Time windows along each edge: $\{(i,j):(i,j)\in\mathcal{A} \setminus \{\delta^+(v_s)\cup\delta^-(v_t)\} \wedge a_i + s_i + \tau_{(i,j)} > b_j\}$
\item Ride-duration limit of each commuter: $\{(i,j),(j,n+i):i\in\mathcal{P} \wedge j\in\mathcal{P}\cup\mathcal{D} \wedge i\neq j \wedge \tau_{(i,j)} + s_j + \tau_{(j,n+i)} > L_i\}$
\item Time windows and ride-duration limits of pairs of trips:
\begin{itemize}
	\item $\{(i,n+j):i,j\in\mathcal{P} \wedge i\neq j \wedge \neg feasible(j\rightarrow i\rightarrow n+j\rightarrow n+i)\}$
	\item $\{(n+i,j):i,j\in\mathcal{P} \wedge i\neq j \wedge \neg feasible(i\rightarrow n+i\rightarrow j\rightarrow n+j)\}$
	\item $\{(i,j):i,j\in\mathcal{P} \wedge i\neq j \wedge \neg feasible(i\rightarrow j\rightarrow n+i\rightarrow n+j) \wedge \neg feasible(i\rightarrow j\rightarrow n+j\rightarrow n+i)\}$
	\item $\{(n+i,n+j):i,j\in\mathcal{P} \wedge i\neq j \wedge \neg feasible(i\rightarrow j\rightarrow n+i\rightarrow n+j) \wedge \neg feasible(j\rightarrow i\rightarrow n+i\rightarrow n+j)\}$
\end{itemize}
\end{enumerate}

\begin{figure}[!t]
	\centering
	\includegraphics[width=0.5\linewidth]{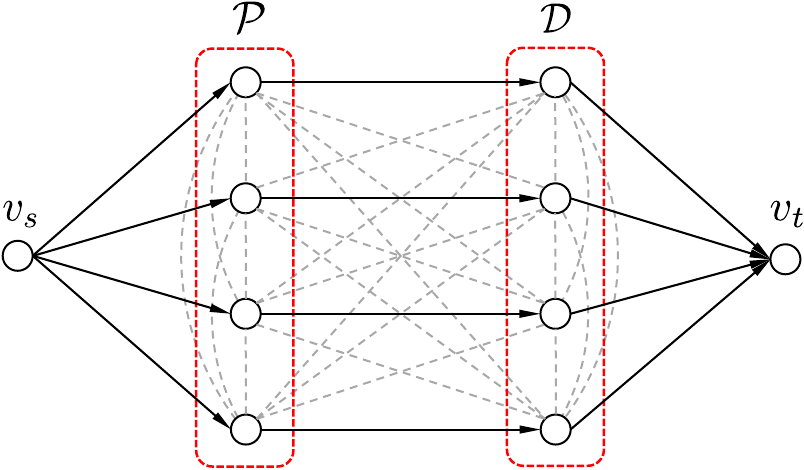}
	\caption{Graph $\mathcal{G}$ (Each Dotted Line Represents a Pair of Bidirectional Edges).}
	\label{fig:graph1}
\end{figure}

\noindent
Note that the sets of edges in (e) utilize the $feasible$ function to
determine if a partial route satisfies time-window and ride-duration
limit constraints. For instance, the first condition indicates that edge
$(i,n+j)$ is infeasible if route $j\rightarrow i\rightarrow
n+j\rightarrow n+i$ is infeasible. Figure \ref{fig:graph1} illustrates
an example of graph $\mathcal{G}$ resulting from the removal of the
infeasible edges.

\paragraph{Big-M Constants} The RMP\textsubscript{CTSPAV} utilizes big-$M$ constants in
constraints \eqref{eqn:ctspav_travel_time_1} and
\eqref{eqn:ctspav_travel_time_2} to enforce the underlying constraints
\emph{only} on selected edges. To ensure that the constants are large
enough to accomplish this goal while not being excessively large so as
to introduce numerical issues, they are defined as follows:
\begin{flalign}
& M_{(i,j)} = \max \{0, b_i + s_i + \tau_{(i,j)} - a_j\} \qquad \forall i,j\in\mathcal{P}\cup\mathcal{D} \\
& \bar{M}_{(i,j)} = \max \{0, b_j - a_i - s_i - \tau_{(i,j)}\} \qquad \forall i\in\mathcal{P}\cup\mathcal{D}, \forall j\in\mathcal{D}
\end{flalign}

\noindent
When the lexicographic objective is considered, edge costs defined in
\eqref{eqn:edge_costs} uses $\hat{\varsigma}_\text{max}$ which is a
constant representing the length of the longest AV route. Since
enumeration of all feasible AV routes in $\mathcal{R}$ is impractical,
a conservative overestimate is used for $\hat{\varsigma}_\text{max}$
to accomplish the lexicographic ordering of the sub-objectives.

\paragraph{Lower Bound} Column-generation procedures are known to have a tailing-off
effect, whereby the rate-of-change of the RMP\textsubscript{CTSPAV}
objective value $z_\text{RMP\textsubscript{CTSPAV}}$ progressively
decreases as $z_\text{RMP\textsubscript{CTSPAV}}$ approaches $z^*$
\citep{lubbecke2005}. To mitigate this effect, a dual lower bound to
$z^*$, $z_\text{LB}$, is defined using the generalized version of the
Lasdon bound \citep{lasdon1970}, i.e.,
\begin{equation} \label{eqn:LasdonLB}
z_\text{LB} = \kappa \bar{c}^*_r
\end{equation}
where $\kappa$ is an upper bound to the number of selected mini routes
in MP\textsubscript{CTSPAV}, $\kappa\geq\sum_{r\in\Omega} X_r$, and
$\bar{c}^*_r$ is the smallest mini-route reduced cost discovered from
PSP\textsubscript{CTSPAV}. For this problem, it is sufficient to take
$\kappa = 2n$. Since the edge costs are all integral, the optimal
objective value of MP\textsubscript{CTSPAV} must also be integral, and
therefore the column-generation iterations can be terminated when
$\ceil{z_\text{RMP\textsubscript{CTSPAV}}} - z_\text{LB} < 1$.

\paragraph{Solving the Subproblem} The label-setting algorithm of \citet{gschwind2015} that is used to solved the
pricing subproblem produces an intermediate set of non-dominated,
feasible mini routes, $\hat{\Omega}_i$ for each graph
$\mathcal{G}_i^+(i\in\mathcal{P}^+)$ and
$\mathcal{G}_i^-(i\in\mathcal{P}^-)$. Instead of considering only the
least-cost route from $\hat{\Omega}_i$, all routes from
$\hat{\Omega}_i$ with negative reduced costs are selected and
introduced into $\Omega'$ to further accelerate the column-generation
convergence. Moreover, since the mini-route search procedure on all graphs are
independent, they are solved concurrently in our
implementation. Finally, $\Omega'$ is initialized with the set of all
direct-trip routes, i.e., it is initialized with $\{i\rightarrow n+i :
i\in\mathcal{P}\}$.

\section{The DARP Column-Generation Procedure}
\label{sec:darp}

The CTSPAV can be viewed as a specialization of the DARP: 
It can be converted into a DARP simply by setting the
time-window of each vehicle at the depot to $\pm \infty$. This section
describes a column-generation procedure, referred to as the DARP
procedure, derived from the algorithm for solving the DARP by
\cite{gschwind2015}. It is the algorithm to which the CTSPAV procedure 
is compared in the computational results
section. At a high level, the DARP procedure is similar to the CTSPAV
procedure as they both use column generation. However, the DARP
procedure fundamentally differs from the CTSPAV procedure in that it adopts the classical column-generation procedure. More specifically, the DARP procedure uses a set-covering restricted master problem RMP\textsubscript{DARP} that only selects AV routes from a set $\mathcal{R}'$ to ensure every trip is covered in the solution. Columns of RMP\textsubscript{DARP} represent AV routes whereas those of RMP\textsubscript{CTSPAV} represent mini routes. The procedure also uses a pricing subproblem PSP\textsubscript{DARP} that searches for feasible AV routes to augment $\mathcal{R}'$. Upon convergence of the column-generation process, the RMP\textsubscript{DARP} is solved as a MIP to obtain an integer solution.\footnotemark

\footnotetext{The branch-and-price approach proposed by
  \cite{gschwind2015} is not considered because it is found to be too
  expensive for the problem instances used in this work. Even the root
  node of the branch-and-price tree cannot be solved within the
  allocated time budget for the real instances considered in this
  paper.}

\subsection{The Master Problem}

The master problem MP\textsubscript{DARP} is a set-covering
formulation that seeks the optimal routing plan for the CTSPAV. It is
defined on the set of all feasible AV routes $\mathcal{R}$ and uses a
binary variable $X_\rho$ that indicates whether route
$\rho\in\mathcal{R}$ is used in the plan. The model is listed in
\eqref{eqn:darp_obj}--\eqref{eqn:darp_route_var}.
\begin{flalign}
& \min \sum_{\rho\in\mathcal{R}} c_\rho X_\rho \label{eqn:darp_obj} \\
& \text{subject to} \nonumber \\
& \sum_{\rho\in\mathcal{R}} a_{i,\rho} X_\rho \geq 1 \qquad \forall i\in\mathcal{P} \label{eqn:darp_route_cover} \\
& X_\rho \in\{0,1\} \qquad \forall \rho\in\mathcal{R} \label{eqn:darp_route_var} 
\end{flalign}
The objective function \eqref{eqn:darp_obj} minimizes the total cost
of the selected routes.  Constant $a_{i,\rho}$ in constraints
\eqref{eqn:darp_route_cover} represents the number of times node $i$ is visited by
route $\rho$. These constraints ensure that each pickup node is
covered in the optimal plan. Our experimental evaluations indicated that a
set-covering formulation produces stronger integer solutions than a
set-partitioning formulation.

To find a routing plan that minimizes vehicle count, the cost $c_\rho$
of each route is set to 1. On the other hand, to find a plan
minimizing total travel distance, the cost $c_\rho$ is set to the
total distance of $\rho$, i.e., $c_\rho =
\sum_{(i,j)\in\rho}\sigma_{(i,j)}$. Finally, to implement a
lexicographic objective that first minimizes vehicle count and then
their total distance, the model is solved twice. The model is first solved to
produce the optimal vehicle count $\chi^*_\text{MIP}$. The constraint
\begin{equation}
\sum_{\rho\in\mathcal{R}} X_\rho = \chi^*_\text{MIP} \label{eqn:darp_vehicle_count} 
\end{equation} 
is then introduced to the model to fix its vehicle count to its
optimal value. The model is then solved again to optimize the
secondary objective. 

While a blended approach similar to that used in the CTSPAV procedure could have also been used here to implement the lexicographic objective, initial experimental evaluations revealed that the greater complexity of PSP\textsubscript{DARP}, which is significantly more expensive than PSP\textsubscript{CTSPAV}, combined with the use of the Lasdon bound \eqref{eqn:LasdonLB} results in a column-generation phase that converges significantly slower. The proposed multi-objective approach, which first just minimizes the vehicle count and therefore uses identical costs for the routes of RMP\textsubscript{DARP} (unlike route costs for the blended approach), permits the use of the dual bound proposed by \cite{farley1990} in the column-generation phase which is  stronger than the Lasdon bound in this setting. This stronger dual bound consequently allows the column-generation termination criterion to be satisfied earlier, thus resulting in a faster converging column-generation phase for the primary objective. And while a similar approach could have also been used for the CTSPAV procedure, the less expensive nature of PSP\textsubscript{CTSPAV} makes a strong dual bound less critical for its column-generation phase which already converges quickly. In the end, the blended and the multi-objective approaches are different yet valid alternatives for implementing the lexicographic objective. The latter, which applies the lexicographic ordering directly, is preferred for the DARP procedure simply because it allows the column-generation phase for the primary objective to converge more quickly in practice and is seen as a necessity to counteract the increased complexity of its pricing subproblem. 

\subsection{The Pricing Subproblem}

The pricing subproblem PSP\textsubscript{DARP} searches for AV
routes with negative reduced costs. Let $\{\alpha_i :
i\in\mathcal{P}\}$ denote the set of optimal duals of constraints
\eqref{eqn:darp_route_cover} and $\beta$ be that of constraint
\eqref{eqn:darp_vehicle_count}. When RMP\textsubscript{DARP} has the
vehicle count-minimization objective, the reduced cost of route $\rho$
is given by:
\begin{equation} \label{eqn:darp_reduced_cost_1}
\bar{c}_\rho = 1 - \sum_{i\in\mathcal{P}} a_{i,\rho} \alpha_i
\end{equation}
When the distance-minimization objective is applied, the reduced cost of $\rho$ is given by:
\begin{equation} \label{eqn:darp_reduced_cost_2}
\bar{c}_\rho = \sum_{(i,j)\in\rho} \sigma_{(i,j)} - \sum_{i\in\mathcal{P}} a_{i,\rho} \alpha_i
\end{equation}
Finally, when constraint \eqref{eqn:darp_vehicle_count} is also present in RMP\textsubscript{DARP} with the distance-minimization objective, the reduced cost of $\rho$ is given by:
\begin{equation} \label{eqn:darp_reduced_cost_3}
\bar{c}_\rho = \sum_{(i,j)\in\rho} \sigma_{(i,j)} - \sum_{i\in\mathcal{P}} a_{i,\rho} \alpha_i - \beta
\end{equation}

The pricing subproblem searches for routes with negative reduced costs
using a graph $\mathcal{G}$ identical to the one described in
Section \ref{sec:CTSPAV}. A reduced cost $\bar{c}_{(i,j)}$ is assigned
to each edge $(i,j)\in\mathcal{A}$, and it is defined according to the
objective function used. For the vehicle count-minimization objective,
$\bar{c}_{(i,j)}$ is given by
\begin{equation} \label{eqn:darp_edge_reduced_cost_1}
\bar{c}_{(i,j)} = 
\begin{cases}
1&\qquad \forall (i,j)\in\delta^+(v_s)\\
-\alpha_i&\qquad \forall i\in\mathcal{P}, \forall j\in\mathcal{N}\\
0&\qquad \forall i\in\mathcal{D}, \forall j\in\mathcal{N}
\end{cases}
\end{equation}
When the distance-minimization objective is used, $\bar{c}_{(i,j)}$ is given by
\begin{equation} \label{eqn:darp_edge_reduced_cost_2}
\bar{c}_{(i,j)} = 
\begin{cases}
\sigma_{(i,j)}-\alpha_i&\qquad \forall i\in\mathcal{P}, \forall j\in\mathcal{N}\\
\sigma_{(i,j)}&\qquad \forall i\in\mathcal{D}\cup\{v_s\}, \forall j\in\mathcal{N}
\end{cases}
\end{equation}
The edge reduced costs are defined so that the total cost of any path
in $\mathcal{G}$ from $v_s$ to $v_t$ is equivalent to the reduced cost
of the path defined in \eqref{eqn:darp_reduced_cost_1}--\eqref{eqn:darp_reduced_cost_3}. In
the presence of constraint \eqref{eqn:darp_vehicle_count} for the
distance-minimization objective, \eqref{eqn:darp_edge_reduced_cost_2}
is still used to define edge reduced costs and $-\beta$ is just added
to the final path cost to obtain the reduced cost defined in
\eqref{eqn:darp_reduced_cost_3}.

Once $\mathcal{G}$ has been set up with the proper edge reduced costs,
the PSP\textsubscript{DARP} just searches for the least-cost path from
$v_s$ to $v_t$ that satisfies the time-window, vehicle-capacity,
pairing, precedence, and ride-duration limit constraints. This
least-cost path is then added to $\mathcal{R}'$ if the cost is
negative. The problem is an SPPRC and is solved using the label-setting dynamic
program proposed by \cite{gschwind2015} which utilizes resource
constraints to enforce the route-feasibility constraints.

In the presence of negative cost cycles in $\mathcal{G}$, the
label-setting algorithm may produce non-elementary paths. Whereas an
additional resource may be introduced to the pricing algorithm to
guarantee path elementarity, the SPPRC then becomes an ESPPRC which is
extremely hard to solve. Therefore, the implementation adopts a
strategy from \cite{ropke2006} and \cite{gschwind2015} which simply
relaxes this elementarity requirement. While doing so theoretically
causes the RMP\textsubscript{DARP} to converge to a weaker primal
lower bound as it now admits a larger set of routes $\mathcal{R}''
\supseteq \mathcal{R}'$, both \cite{ropke2009} and \cite{gschwind2015}
have found that the resulting lower bound is only slightly weaker in
practice.

\subsection{Implementation Strategies}

Similar to the CTSPAV procedure, a dual lower bound is maintained
during the column-generation procedure to mitigate its tailing-off
effect. When the vehicle count-minimization objective is active, the
lower bound proposed by \cite{farley1990} is adopted since the cost of
each route is identical. It is given by:
\begin{equation}
z_\text{LB}' = \frac{z_\text{RMP\textsubscript{DARP}}}{1 - \bar{c}_\rho^*}
\end{equation}
where $z_\text{RMP\textsubscript{DARP}}$ is the objective value of
RMP\textsubscript{DARP} at the end of each iteration and
$\bar{c}_\rho^*$ is the smallest route reduced cost discovered by
PSP\textsubscript{DARP}. The column generation is then terminated when
$\ceil{z_\text{RMP\textsubscript{DARP}}} - z_\text{LB}' < 1$.  To
accomplish the lexicographic objective during the column generation,
constraint \eqref{eqn:darp_vehicle_count} is introduced to
the RMP\textsubscript{DARP} after the primary objective has converged with
its right-hand side set to
$\ceil{z_\text{RMP\textsubscript{DARP}}}$. The objective function is
then switched to distance minimization after which the column generation
is resumed. For the distance-minimization objective, the generalized
Lasdon bound defined in \eqref{eqn:LasdonLB} is used as a dual lower
bound, and the column generation is terminated when
$\ceil{z_\text{RMP\textsubscript{DARP}}} - z_\text{LB} < 1$.

Our implementation also incorporates the interior-point, dual-stabilization 
method proposed by \cite{rousseau2007} to accelerate
the column-generation convergence. In addition, all non-dominated
routes with negative reduced costs resulting from the label-setting
algorithm in PSP\textsubscript{DARP} are added to the master problem. The column-generation phase is seeded with routes $\{i\rightarrow n+i\rightarrow 2n+i\rightarrow 3n+i:i\in\mathcal{P}^+\}$. Since this set of routes represents a feasible solution to the problem, it guarantees the existence of a feasible integer solution to the MIP for both the distance-minimization and the primary lexicographic objective. Recall that for the secondary lexicographic objective, the right-hand side of the introduced constraint \eqref{eqn:darp_vehicle_count} is set to $\chi^*_\text{MIP}$ which represents the objective value of the optimal integer solution for the primary objective. The existence of this feasible integer solution guarantees that one also exists for the MIP with the secondary objective.

Non-elementary routes produced in the column-generation phase are
removed prior to solving the RMP\textsubscript{DARP} as a MIP.
Repeated nodes are identified from each non-elementary route and only
the first instance of each repeated node is preserved in the route
(subsequent instances are removed). The resulting route is feasible as
the non-elementary version already satisfies the time-window,
vehicle-capacity, pairing, precedence, and ride-duration limit
constraints. Similarly, since a set-covering formulation is used, a
node may be visited by multiple routes in the integer solution. This
is fixed by simply preserving the node in an arbitrarily selected
route and removing it from the others. Since the travel distances
satisfy the triangle inequality, this step only shortens the affected
routes and hence improves the total travel distance of the solution
while maintaining its vehicle count.

\section{Case Study and Experimental Results}
\label{sec:results}

This section evaluates the potential benefits of autonomous vehicles
on a real case study. It also reports a variety of experimental
results on the efficiency of the optimization algorithms.

\subsection{The Dataset and Construction of Problem Instances}

\begin{figure}[!t]
	\centering
	\includegraphics[width=1.0\linewidth]{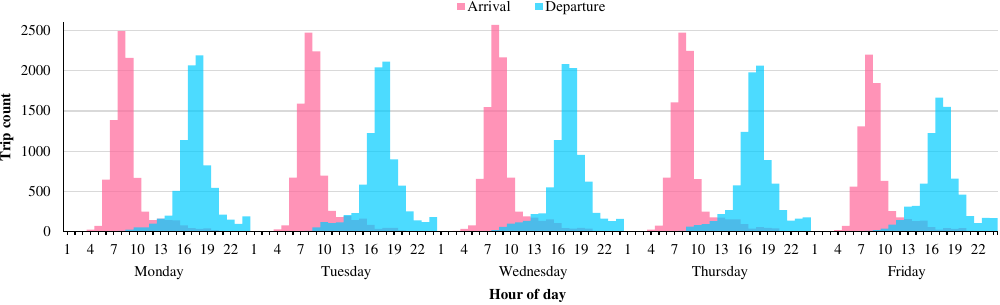}
	\caption{Distribution of Arrival and Departure Times of Commute Trips from Parking Structures over Week 2 of April 2017.}
	\label{fig:arrivedepart_distribution}
\end{figure}

The performance of the CTSPAV and DARP procedures are evaluated on
problem instances derived from the real-world, commute-trip dataset
used by \citet{hasan2018}. The dataset consists of daily arrival and
departure information of a total of 15,000 commuters traveling to 15
parking structures located in downtown Ann Arbor, Michigan. Collected
throughout the month of April 2017, this information includes the
exact arrival and departure times of every commuter to and from each
parking structure. Joining this information with the addresses of the
commuter homes ---which are situated within the city as well as its
surrounding region, an area that spans 13,000 square miles--- yields
detailed information about an average of 9,000 commute trips per
weekday. The results in this section focus on the busiest days
(Monday--Thursday) of the busiest week of the month (week 2). Figure
\ref{fig:arrivedepart_distribution} shows the distribution of arrival
and departure times of these trips to and from the parking structures
over this week. The travel patterns display remarkably consistent
distributions over the different weekdays. The arrival and departure
times peak at 6--9 am and 4--7 pm respectively, highlighting the
typical peak commute hours.

For additional perspectives, the trips are partitioned into two sets:
the approximately 2,200 daily trips made by commuters living within
Ann Arbor city limits (the region bounded by highways US-23, M-14, and
I-94), and the remainder made by commuters living outside the
region. The rationale behind this split is to distinguish between the
results of trips for commuters living nearer to the downtown area,
where the parking structures are located, and those for commuters
living further away. The trips from each set are then further
partitioned into smaller problem instances using the clustering
algorithm from \citet{hasan2020} which groups at most $N$ commuters
together based on the spatial proximity of their home locations. The
clusters are thought of as ``neighborhoods'' within which trip sharing
is done exclusively. This notion of breaking down a problem into
smaller instances is in line with the conclusion by \citet{agatz2012},
that effective decomposition techniques are necessary for the
computational feasibility of large-scale problems.


Several assumptions were made to characterize the nature of the
requests submitted to the trip-sharing platform. First, each rider
$i$, when requesting a trip, specifies a desired arrival time $at^+_i$
to the inbound trip destination and a desired departure time $dt^-_i$
from the outbound trip origin. This assumption is consistent with
those in the literature on the DARP, e.g., \citet{jaw1986},
\citet{cordeau2003}, and \citet{cordeau2006}. Secondly, each rider
tolerates a maximum shift of $\pm\Delta$ to the desired
times. Therefore, if the arrival and departure times at the parking
structure from the dataset are considered as the desired times, then a
delivery time upper bound of $b_i$ = $at^+_i + \Delta$ can be
associated with each $i\in\mathcal{D}^+$ and a pickup time window of
$[a_i,b_i] = [dt^-_i - \Delta,dt^-_i + \Delta]$ can be associated with
each $i\in\mathcal{P}^-$. Consequently, a departure time window of
$[a_i,b_i] = [b_{n+i} - s_i - L_i - 2\Delta,b_{n+i} - s_i - L_i]$ is
associated with each $i\in\mathcal{P}^+$ and a delivery time upper
bound of $b_i = b_{i-n} + s_{i-n} + L_{i-n}$ is associated with each
$i\in\mathcal{D}^-$. Finally, each rider $i$ will tolerate at most an
$R\%$ extension to her direct-trip ride duration, i.e., $L_i =
(1+R)\cdot \tau_{(i,n+i)}$ for each
$i\in\mathcal{P}$. \citet{hunsaker2002} also used a similar assumption. This use of a single multiplication factor for modeling the maximum tolerable ride-duration (versus a more fine-grained approach) is seen as a practical necessity to cater to the massive volume of trips considered in the case study, and it is seen as a more realistic alternative to the approach used in other works on the DARP (e.g., \cite{cordeau2006} used an identical duration for every trip). Results of sensitivity analyses on both $\Delta$ and $R$ are provided in later sections to demonstrate their effect on the final results. 

\paragraph{Depot Configurations}

This paper explores two hypothetical depot configurations: (1) a
central depot configuration in which all neighborhoods are served by
vehicles from a single, centralized depot, and (2) a local depot
configuration whereby each neighborhood is served by a local depot
situated within the neighborhood itself. For the first scenario, the
largest parking structure from the dataset considered is arbitrarily
designated as the central depot. In the second scenario, the home
address $l_c$ that is nearest to every other location within cluster
$\mathcal{P}_c$ is selected as the hypothetical local depot location,
i.e.,
\begin{equation}
l_c = \argmin_{i\in\mathcal{P}_c} \bigg\{\sum_{j\in\mathcal{P}_c\setminus\{i\}} \sigma_{(i,j)} + \sigma_{(j,i)}\bigg\}.
\end{equation}  

\subsection{Experimental Setup and Parameters}

The GPS coordinates of every address considered are geocoded using
Geocodio, while the shortest path, travel time, and travel distance
between any two locations are estimated using GraphHopper's Directions
API that uses OpenStreetMap data. All algorithms are implemented in
C++ with parallelization being handled with OpenMP. The label-setting
algorithm of \citet{gschwind2015} is implemented using the
resource-constrained shortest path framework from the Boost Graph
Library (version 1.70.0), while all LPs and MIPs are solved with
Gurobi 8.1.1. Every problem instance is solved on a compute cluster,
utilizing 12 cores of a 3.0 GHz Intel Xeon Gold 6154 processor and 32
GB of RAM. A total time budget of 2 hours is allocated for each
instance; 1 hour for the column-generation phase and another 1 hour
for solving the MIP for the CTSPAV procedure. The same total time
budget is allocated for the DARP procedure. However, since initial
evaluations showed that it requires more time for the
column-generation phase and less for solving the MIP, 1.5 hours is
allocated for its column-generation phase and 0.5 hours for solving
its MIP. All reported results are from the best feasible solution
obtained within the time limit.

The experiments consider the use of autonomous cars and hence
$K=4$. Initial empirical evaluations found that $N=100$ for the
clustering algorithm produces neighborhoods that are sufficiently
large to provide ample opportunities for trip sharing while not
producing intractable problem instances. This setting was thus used to
generate the problem instances for all experiments. Finally, unless
otherwise stated, values of $\Delta = 10$ mins and $R = 50\%$ are used
in all experiments. \ref{sec:gaps} summarizes results on the
optimality gaps and computation times of the procedures. Both the
appendix and the next subsection summarizes the results of all problem
instances obtained from applying the clustering algorithm on all commute
trips from the Wednesday of week 2 (which generated 22 and 68 clusters inside
and outside the city respectively).

\subsection{Performance Comparison of the CTSPAV and DARP Procedures}

\begin{figure}[!t]
	\centering
	\includegraphics[width=0.75\linewidth]{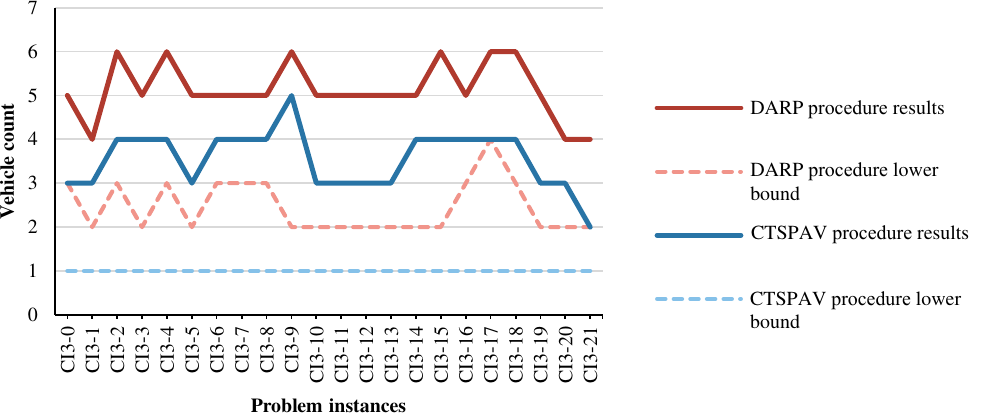}
	\caption{Comparison of Vehicle Count Results for Problem
          Instances Inside City Limits with a Lexicographic Objective
          and a Central Depot Configuration.}
	\label{fig:vc_darp_ctspav_inside}
\end{figure}

\begin{figure}[!t]
	\centering
	\includegraphics[width=1.0\linewidth]{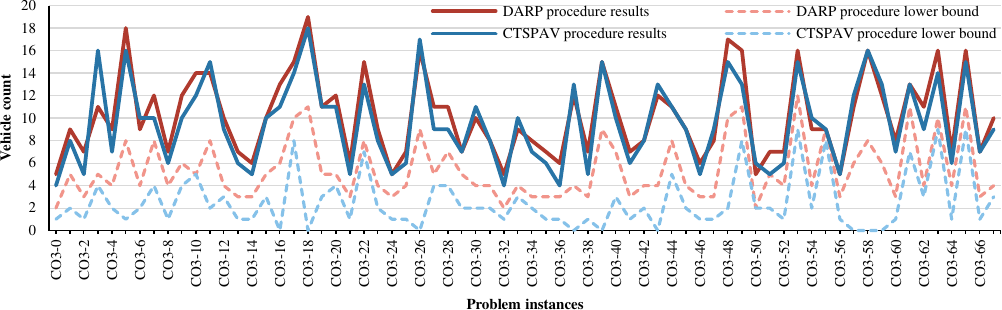}
	\caption{Comparison of Vehicle Count Results for Problem
          Instances Outside City Limits with a Lexicographic Objective
          and a Central Depot Configuration.}
	\label{fig:vc_darp_ctspav_outside}
\end{figure}

This section presents a comparison of the results of the CTSPAV
procedure against those of the DARP, particularly the final objective
values from their MIPs together with the corresponding lower bounds of
each procedure for every problem instance considered. Figure
\ref{fig:vc_darp_ctspav_inside} compares their vehicle count results
when the lexicographic objective and central depot configuration are
used on clusters inside the city, whereas Figure
\ref{fig:vc_darp_ctspav_outside} does the same for clusters
outside. The horizontal axis lists the problem instances while the
vertical displays the vehicle counts. Figure
\ref{fig:vc_darp_ctspav_inside} shows the typical trend: the DARP
procedure produces stronger lower bounds and the CTSPAV procedure
produces smaller vehicle counts. This trend carries over to Figure
\ref{fig:vc_darp_ctspav_outside}. Even though there are a few
instances where the DARP procedure outperform the CTSPAV procedure,
overall the CTSPAV matches or outperforms the DARP on more than 80\%
of the instances shown in the figure. Vehicle count results of
instances with the local depot configuration exhibit similar trends
and are not shown here.

Table \ref{tab:my-table} then compares the aggregated vehicle count results from all clusters for the same lexicographic objective and central depot configuration. This comparison is especially important as the results of all clusters are summarized through aggregation in the case study. The table shows a significant difference between the results of the two procedures, with the DARP procedure producing aggregated vehicle counts that are 45\% larger inside the city, 6\% larger outside the city, and 11\% larger overall. This significant difference makes another strong case for using results of the CTSPAV procedure in the case study.

\begin{table}[!t]
	\centering
	\caption{Difference in Aggregate Vehicle Counts of CTSPAV and DARP Procedures}
	\label{tab:my-table}
		\begin{tabular}{cccc}
			\hline \noalign{\smallskip}
			\multirow{2}{*}{Location} & \multicolumn{3}{c}{Aggregate vehicle count} \\
			\cline{2-4} \rule{0pt}{12pt}
			& CTSPAV procedure & DARP procedure & Percentage difference \\
			\hline \noalign{\smallskip}
			Inside & 78 & 113 & +45\% \\
			Outside & 652 & 694 & +6\% \\
			Combined & 730 & 807 & +11\% \\
			\hline
		\end{tabular}
\end{table}

\begin{figure}[!t]
	\centering
	\includegraphics[width=0.75\linewidth]{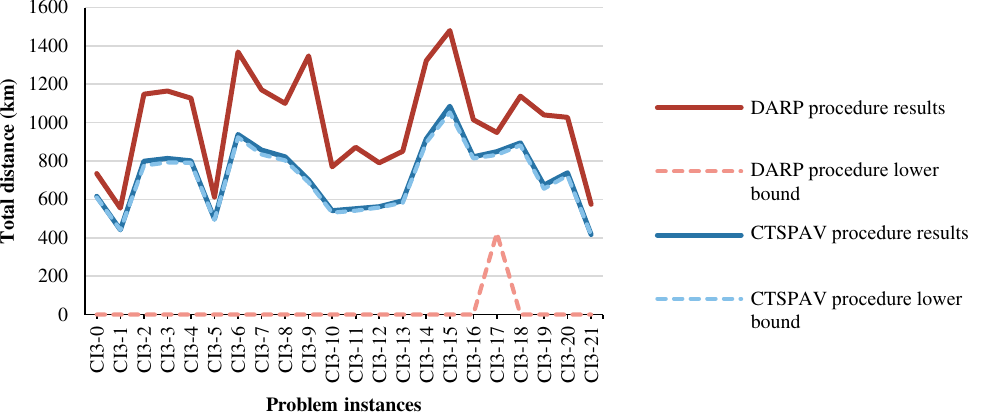}
	\caption{Comparison of Procedure Total Distance Results for Problem Instances Inside City Limits with Distance-Minimization Objective and Central Depot Configuration.}
	\label{fig:dist_darp_ctspav_inside}
\end{figure}

\begin{figure}[!t]
	\centering
	\includegraphics[width=1.0\linewidth]{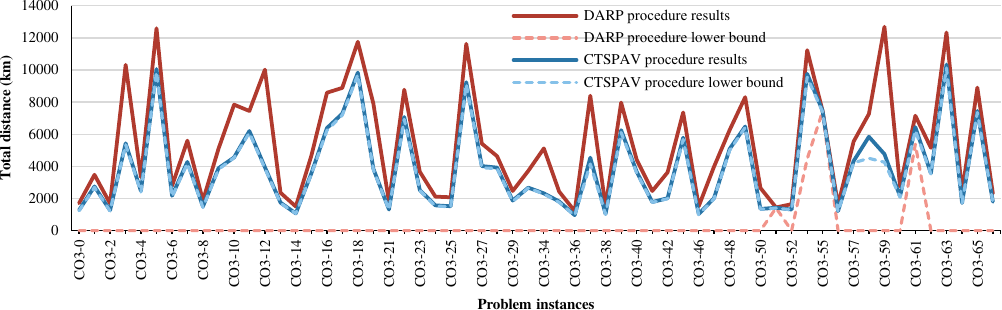}
	\caption{Comparison of Procedure Total Distance Results for Problem Instances Outside City Limits with Distance-Minimization Objective and Central Depot Configuration.}
	\label{fig:dist_darp_ctspav_outside}
\end{figure}

Figures \ref{fig:dist_darp_ctspav_inside} and
\ref{fig:dist_darp_ctspav_outside} compare the results when minimizing
total distance under a central depot configuration for clusters inside
and outside the city respectively. The horizontal axis lists the
problem instances while the vertical displays the total distances. For
a few clusters outside the city, even a single column-generation
iteration of the DARP procedure cannot be completed within the time
limit. The comparison for these instances are therefore excluded from
Figure \ref{fig:dist_darp_ctspav_outside}. Similar to Figures
\ref{fig:vc_darp_ctspav_inside} and \ref{fig:vc_darp_ctspav_outside},
the CTSPAV procedure produces stronger total distance results. In
fact, it outperforms the DARP procedure in all problem instances
considered. In addition, the procedure also appears to produce
stronger lower bounds for this objective function. There is a caveat
to this observation however. The weak lower bounds of the DARP
procedure in this case can be attributed mainly to its
column-generation phase not converging within its (longer) time
limit. When column-generation does converge, e.g., in instances CO3-51
and CO3-55 in Figure \ref{fig:dist_darp_ctspav_outside}, the lower
bounds produced are comparable to those of the CTSPAV procedure. The
figures also highlight the excellent optimality gap of the CTSPAV
procedure for this objective function. Once again, total distance
results of instances with the local depot configuration are not
summarized here as they display similar trends. Since the CTSPAV
procedure consistently produces the stronger final results, it is used
to obtain the results for the subsequent case study.

\subsection{Vehicle and Travel Distance Reduction Results}

\begin{figure}[!t]
	\centering
	\begin{minipage}{.49\textwidth}
		\centering
		\includegraphics[width=1.0\linewidth]{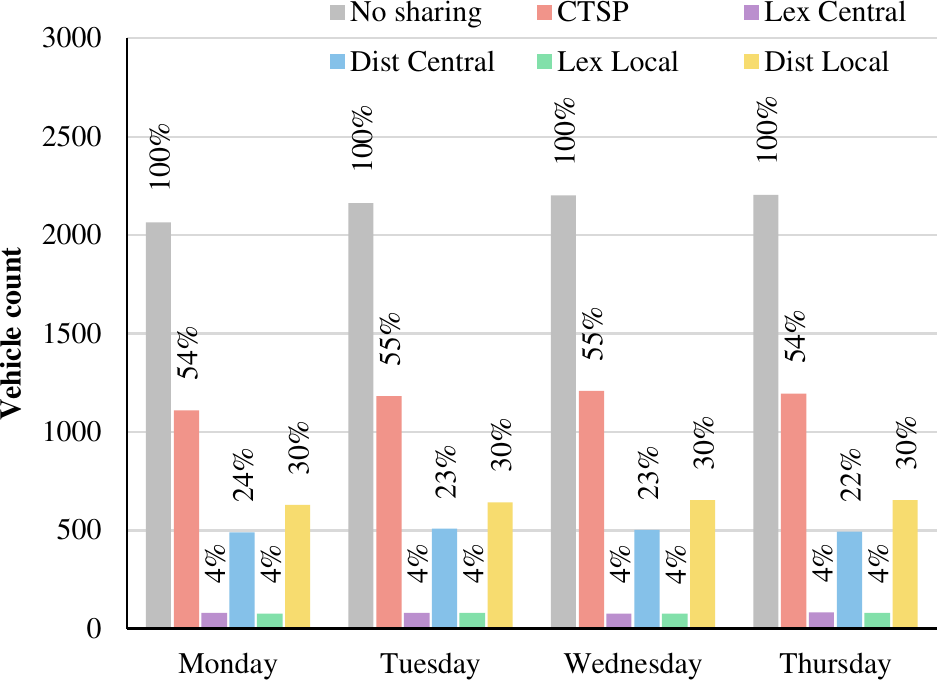}
		\caption{Aggregate Vehicle Count Results from All Clusters Inside City Limits.}
		\label{fig:vehcount_inside}
	\end{minipage}
	\hspace{0.0\textwidth}
	\begin{minipage}{.49\textwidth}
		\centering
		\includegraphics[width=1.0\linewidth]{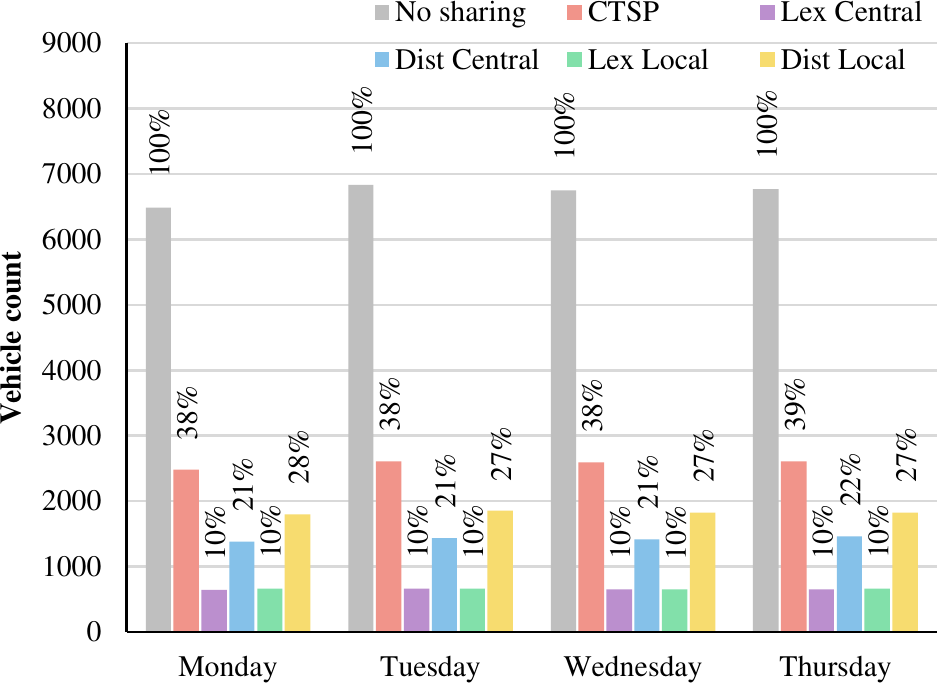}
		\caption{Aggregate Vehicle Count Results from All Clusters Outside City Limits.}
		\label{fig:vehcount_outside}
	\end{minipage}
\end{figure}

Figures \ref{fig:vehcount_inside} and \ref{fig:vehcount_outside} show
aggregated vehicle count (VC) results of all clusters inside and
outside city limits respectively for the first four weekdays of week
2. Each figure shows VC results for every combination of objective
function (lexicographic or distance minimization) and depot
configuration (central or local) for the CTSPAV: they are labeled ``Lex
Central'', ``Dist Central'', ``Lex Local'', and ``Dist Local'' respectively. Each
figure also shows VC results of trips under no-sharing conditions and of the
CTSP described by \citet{hasan2020} for additional perspective. The
percentages of the VCs for each method as a fraction of the no-sharing
VC are also included. The figures indicate that the CTSPAV
consistently requires fewer vehicles than the CTSP to cover all trips
regardless of the location of the clusters, objective function, and
depot configuration. This is not surprising, as the CTSPAV addresses
the key limitation of the CTSP. As mentioned by \citet{hasan2020},
routes of the CTSP are relatively short as the number of locations
they can visit are limited by the ride-duration constraint of their drivers,
as well as the time windows at the origins and destinations of the
driver trips. In contrast, AVs are not subject to these
limitations. Therefore, they can travel back and forth between the
parking structures and the neighborhoods to serve trips throughout the
day, consequently allowing fewer vehicles to be utilized to cover the
same amount of trips. {\em What is striking however is the magnitude
  of the reduction in the number of vehicles: the VCs are reduced by 96\% 
  and 90\% inside and outside city limits respectively when using AVs.}

\begin{figure}[!t]
		\centering
		\includegraphics[width=0.5\linewidth]{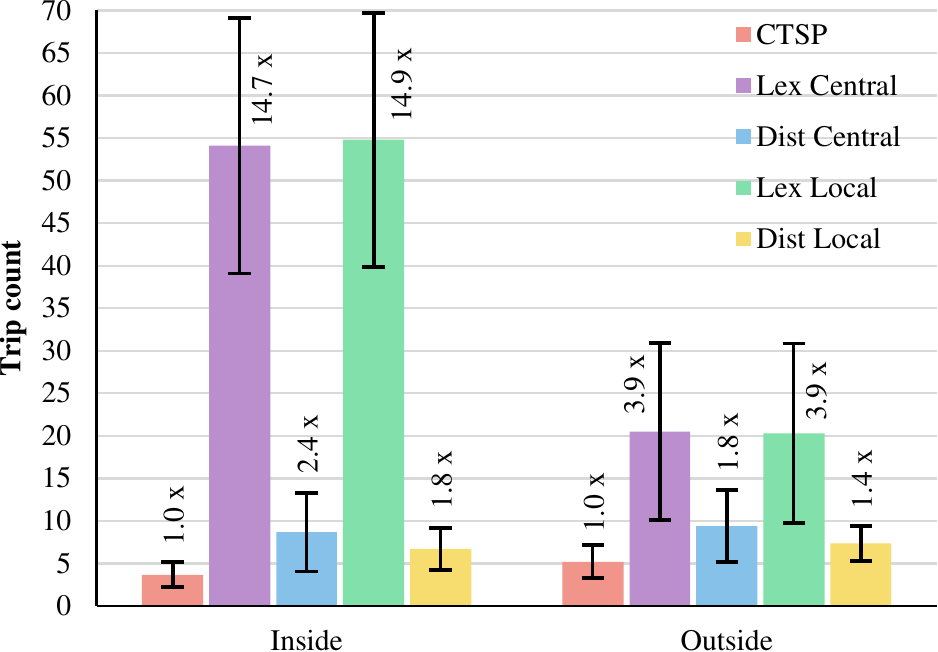}
		\caption{Average Number of Trips Served by Routes of Each Method.}
		\label{fig:route_tripcount}
\end{figure}

Results for the different objective functions show that the vehicle
reduction of the lexicographic objective is significantly better than
that of the distance minimization. Again, this result is not surprising
as reducing VC is the primary goal of the former objective function,
whereas it is not a consideration in the latter.  The discrepancy between the vehicle reductions
inside and outside city limits can be attributed to the larger
distance between the parking structures and the neighborhoods outside
city limits. As a result, the AVs serving these neighborhoods spend a
larger fraction of their time in the transit phase, therefore limiting
the number of trips they could serve in a day. This factor is further
highlighted in Figure \ref{fig:route_tripcount} which summarizes the
average number of trips served by the routes of each method. The
figure also shows each count as a multiplicative factor of the count
for the CTSP, and the error bars depict the standard deviations of the
trip counts. For the CTSPAV with the lexicographic objective, the
routes from the clusters outside the city visit significantly fewer
nodes than those inside the city as their transit phases are
longer. However, regardless of the position of the clusters, routes of
the CTSPAV consistently cover more trips on average than those of the
CTSP. {\em In fact, inside city limits, routes of the CTSPAV with the
lexicographic objective serve, on average, an order of magnitude more
trips than those of the CTSP.}

\begin{figure}[!t]
	\centering
	\begin{minipage}{.49\textwidth}
		\centering
		\includegraphics[width=1.0\linewidth]{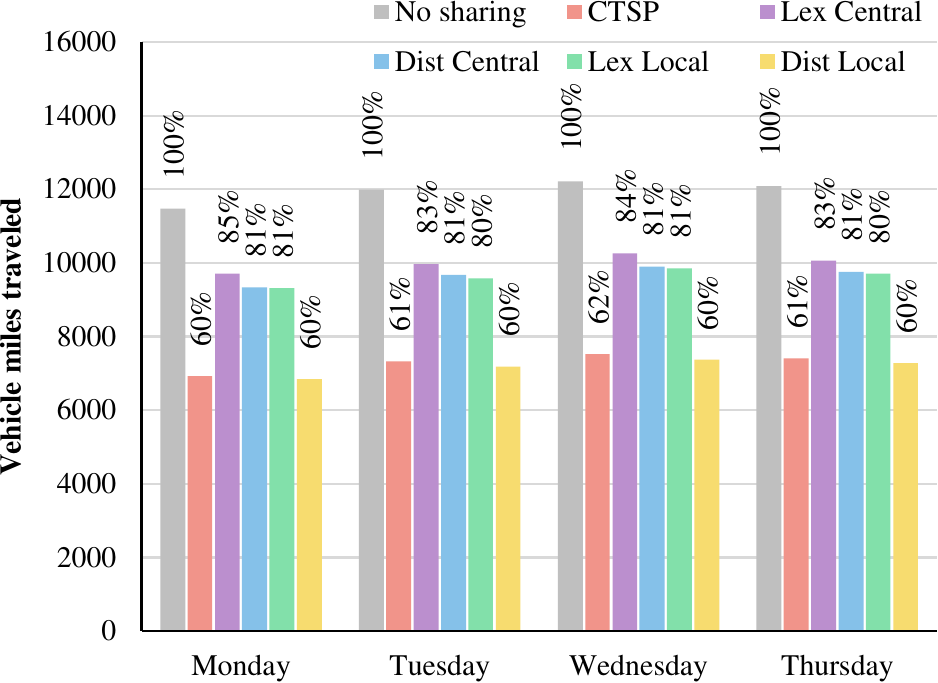}
		\caption{Aggregate Vehicle Miles Traveled from All Clusters Inside City Limits.}
		\label{fig:dist_inside}
	\end{minipage}
	\hspace{0.0\textwidth}
	\begin{minipage}{.49\textwidth}
		\centering
		\includegraphics[width=1.0\linewidth]{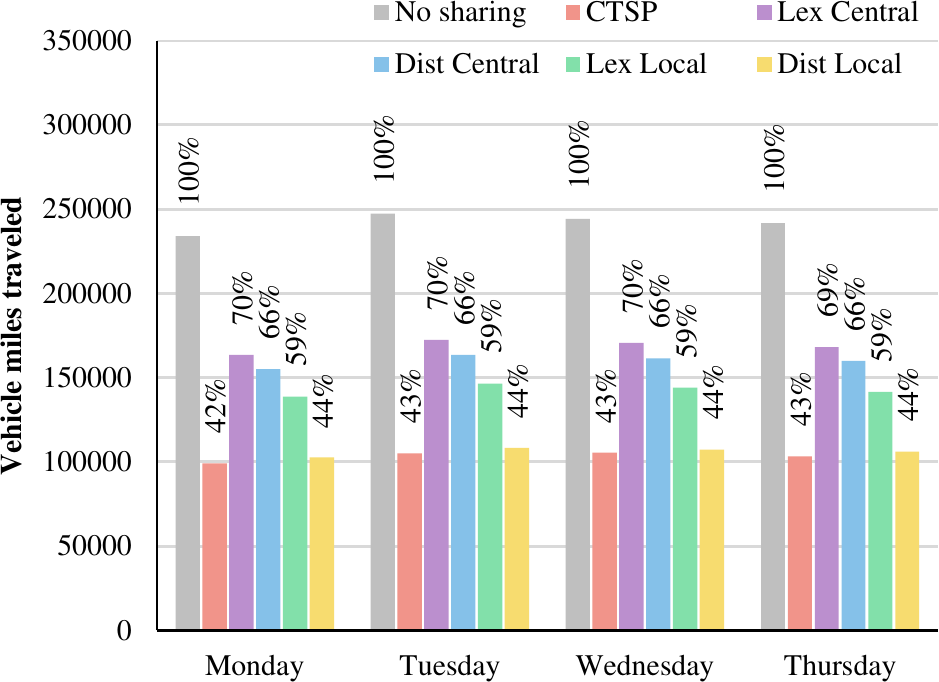}
		\caption{Aggregate Vehicle Miles Traveled from All Clusters Outside City Limits.}
		\label{fig:dist_outside}
	\end{minipage}
\end{figure}

Figures \ref{fig:dist_inside} and \ref{fig:dist_outside} summarize the
corresponding aggregated vehicle miles traveled (VMT) for all clusters
inside and outside city limits respectively. Similar to Figures
\ref{fig:vehcount_inside} and \ref{fig:vehcount_outside}, they show
CTSPAV results for every combination of objective function and depot
configuration, as well as results of the CTSP and of trips under
no-sharing conditions for additional perspectives. Similarly, the
percentages represent the VMT of each method as a fraction of the
no-sharing VMT. The VMT percentage of each method outside the city limits
is consistently smaller than those inside. This can be attributed to
the neighborhoods outside the city being further away from the parking
structures.

When a central depot is used for the CTSPAV, VMT is reduced by
15--19\% inside the city (resp., 30--34\% outside the city). Although
sizeable, the reduction is not as significant as that of the CTSP, as
the AVs have to travel back and forth between the neighborhoods and
the parking structures which increases their total travel
distance. When the different objective functions for the central depot
configuration are compared, the distance-minimization objective only
improves VMT by 2--4\% at the expense of significantly larger vehicle
counts. VMT is further reduced when local depots are used, as each AV
has to travel a shorter distance from the depot to reach its first
pickup location and from its last drop-off location back to the depot
since both these locations are in the neighborhoods. When the
different objective functions for this depot configuration are
compared, it can be seen that the distance-minimization objective
significantly improves VMT by 15--21\%, once again at the expense of a
significant increase in vehicle count. In fact, the VMT of the
distance-minimization objective is comparable to that of the
CTSP. This can be explained by referring to Figure
\ref{fig:route_tripcount}, which shows their trip counts to be similar.
In other words, the VMT of the CTSPAV with the distance-minimization
objective and local depot configuration is comparable to the CTSP
because its routes are serving fewer trips per day.

\begin{figure}[!t]
		\centering
		\includegraphics[width=0.5\linewidth]{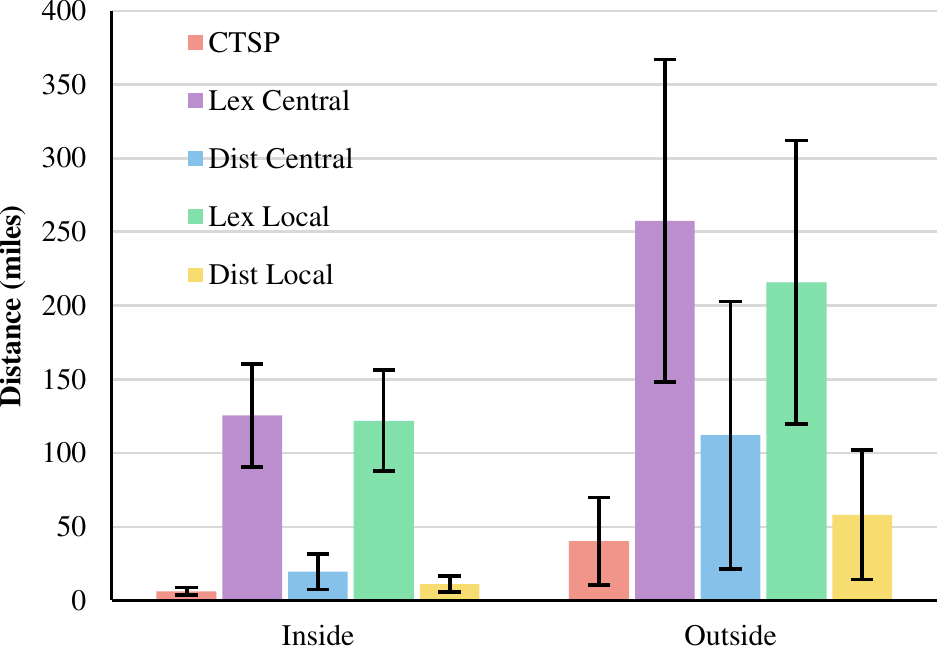}
		\caption{Average Distance Traveled Per Vehicle of Each Method.}
		\label{fig:route_distance}
\end{figure}

Figure \ref{fig:route_distance} provides further insight into the
daily average of the distance traveled per vehicle for each
method. The error bars depict the standard deviations for each
method. As expected, average travel distances are larger for clusters
outside the city for each method. They are also larger for the CTSPAV
when the lexicographic objective is used as the AVs travel back and
forth more between the neighborhoods and parking structures in order
to reduce vehicle count.

\begin{figure}[!t]
	\centering
	\includegraphics[width=1.0\linewidth]{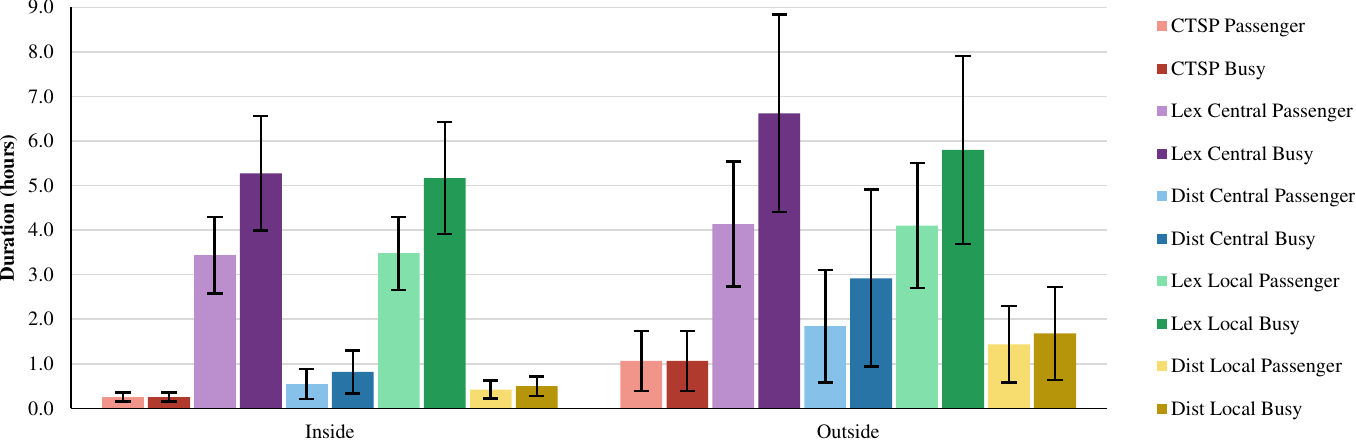}
	\caption{Average Passenger and Busy Ride Durations of the Routes of Each Method.}
	\label{fig:route_passenger_busyduration}
\end{figure}

Figure \ref{fig:route_passenger_busyduration} takes a look at the
daily average passenger and busy ride durations. The error bars again represent the standard deviations for each
method. The passenger duration is defined as the total
duration of the day during which at least one passenger is on the
vehicle, whereas the busy duration is the passenger duration of a
vehicle combined with the duration spent traveling between locations
with no passengers. The complement of the busy duration, the idle
duration, is therefore the duration of the day during which a vehicle
is parked at a depot/at home/at a parking structure, combined with the
duration spent waiting to pick up passengers while being empty. The
passenger and busy durations of the CTSP are therefore identical as
its vehicles are driven by the commuters themselves. In other words,
they never travel without any passengers. For the CTSPAV, its busy
duration is longer than its passenger duration as the latter is a
subset of the former. In fact, the difference between the two
represent the duration spent by the vehicle traveling without any
passengers. For the lexicographic objective, this duration takes up a
large fraction of its busy time, which supports earlier claims that
vehicles under this configuration do more back and forth traveling to
reduce its vehicle count.

The average durations of each method for clusters outside city limits
are longer than corresponding durations inside as the vehicles must
travel farther between the neighborhoods and the parking
structures. The CTSP also has the shortest passenger and busy
durations, which are consistent with the results from Figure
\ref{fig:route_tripcount} which showed the CTSP having the lowest
average trip count. For the various configurations of the CTSPAV, the
ones utilizing the lexicographic objective produce the largest
passenger and busy durations, also consistent with their average trip
count results from Figure \ref{fig:route_tripcount}. The results also
reveal another drawback of the distance-minimization objective. Not
only does the configuration require relatively larger vehicle counts
to cover the same number of trips, but its vehicles are less busy
throughout the day than those of the lexicographic objective. In other
words, the vehicles spend a longer time every day being idle, which
somewhat defeats the purpose of utilizing AVs in the first place.

\begin{figure}[!t]
	\centering
	\begin{minipage}{.49\textwidth}
		\centering
		\includegraphics[width=1.0\linewidth]{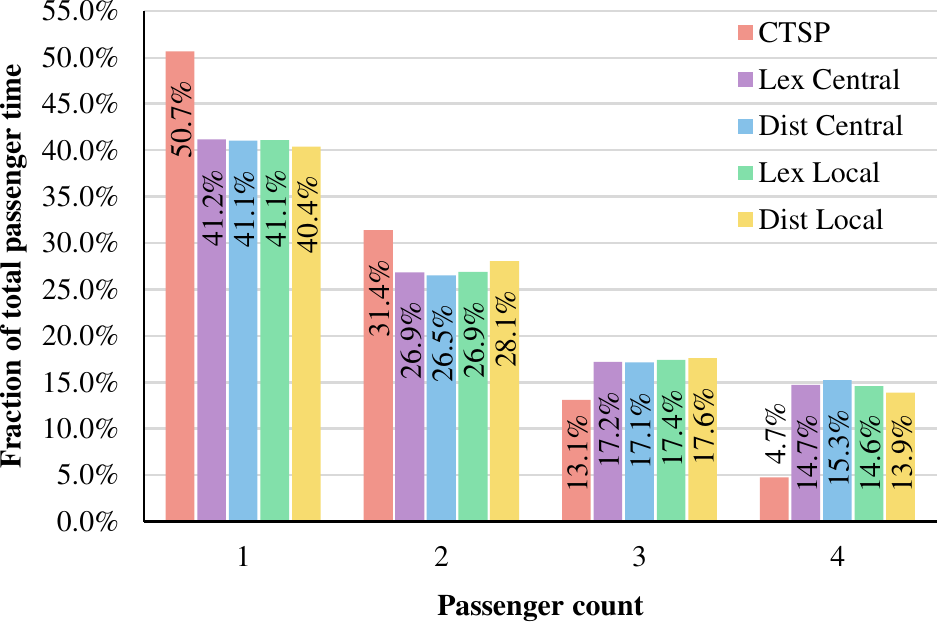}
		\caption{Fraction of Total Passenger Time Spent Serving 1, 2, 3, and 4 Passengers for Clusters Inside City Limits.}
		\label{fig:time_distribution_inside}
	\end{minipage}
	\hspace{0.0\textwidth}
	\begin{minipage}{.49\textwidth}
		\centering
		\includegraphics[width=1.0\linewidth]{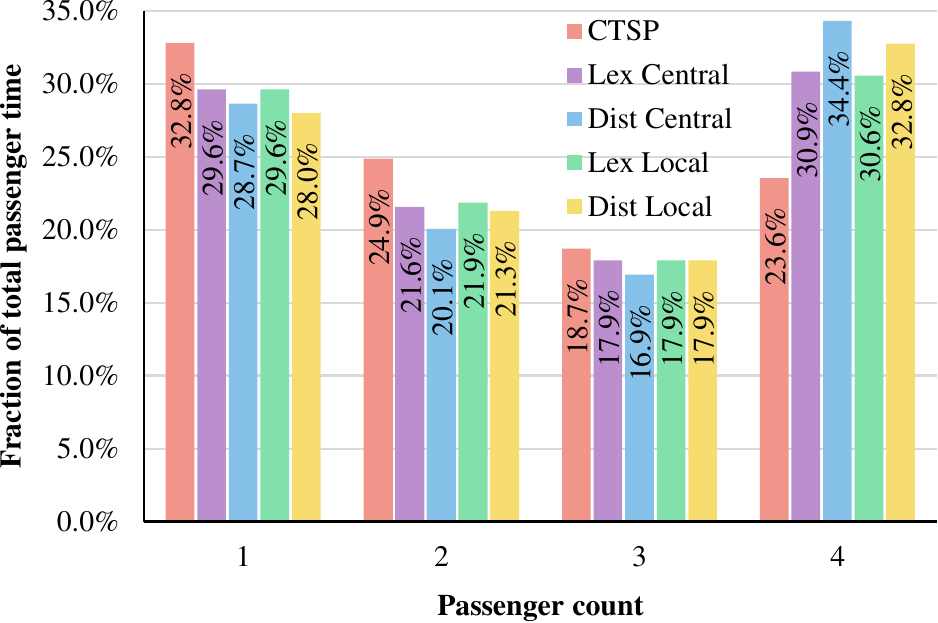}
		\caption{Fraction of Total Passenger Time Spent Serving 1, 2, 3, and 4 Passengers for Clusters Outside City Limits.}
		\label{fig:time_distribution_outside}
	\end{minipage}
\end{figure}

Figures \ref{fig:time_distribution_inside} and
\ref{fig:time_distribution_outside} provide a deeper look into how the
passenger durations are spent by the routes of each method for
clusters inside and outside city limits respectively. It shows the
fraction of the total passenger duration spent serving one, two,
three, or four passengers for each configuration. Inside city limits,
these fractions progressively decrease as the number of passengers
increases. However the CTSPAV spends a larger fraction of its time
serving three or four passengers than the CTSP, again reinforcing the
advantage of the CTSPAV. This is due to two factors. First, mini
routes in the CTSP must start and finish with the same driver. Second,
the CTSP must also synchronize the inbound and outbound routes, since
the same drivers are used for both. Because the schedules of the
passengers very often differ, some of the mini routes cannot be used
by the CTSP but can by the CTSPAV. This observation is carried over to
clusters outside the city, whereby every configuration of the CTSPAV
spends most of its time serving four passengers, while the same cannot
be said for the CTSP. This observation can be attributed the longer
transit phase of the routes for clusters outside the city combined
with the vehicles being used to their full capacity during the transit
phase by the CTSPAV.

In summary, the lexicographic objective for the CTSPAV has the
greatest vehicle reduction potential, and the configuration of the
depot does not appear to affect this potential. Its routes
consistently cover the largest amount of trips on average, and they
also have the longest passenger and busy durations, which are all
desirable characteristics for an effective AV trip-sharing
platform. The local depot configuration for the lexicographic
objective does produce slightly better VMT results. This small
benefit, however, may be outweighed by the logistical benefits of
having a central depot, e.g., the convenience and cost effectiveness
of having a central location for maintaining and refueling/recharging
all AVs. The CTSPAV with the distance minimization objective and local
depot configuration consistently produces the lowest VMT, however as
mentioned earlier, the result is obtained at the expense of higher
vehicle counts. Besides that, as shown in Figures
\ref{fig:route_tripcount} and \ref{fig:route_passenger_busyduration},
the AVs also serve relatively fewer trips and spend more of their times
being idle every day under this configuration.

\subsection{Cost Analysis}

This section reports a simple, coarse analysis to estimate the cost of operating
the trip-sharing platform over a 5- and 10-year period. The analysis is not intended to be a sophisticated or complete measure of the total cost of the platform per se; instead, it is aimed to obtain a rough understanding of the trade-off between vehicle and operating costs of the platform. The analysis focuses only on the CTSPAV with the central depot configuration, as it is preferred over the local depot configuration due to its
aforementioned logistical benefits. The analysis first considers
vehicle-related costs (referred to simply as vehicle costs) and then
operating costs.

The vehicle cost is age-related and consists of the vehicle
depreciation over $y$ years and a distance-related cost. This last
cost consists of an average fuel cost of \$0.15 and an average
vehicle-value depreciation of \$0.08 per mile traveled.  An
exponential decay function is used to model how a vehicle value
depreciates with age. More precisely, its depreciation $\gamma$ over
$y$ years is given by:
\begin{equation}
\gamma = p - p(1 - \nu)^y
\end{equation}
where $p$ is the vehicle's initial price and $\nu$ is its annual
depreciation rate. This analysis uses a depreciation rate of 24\% for
the first year and 15\% for subsequent years. The total vehicle cost
over $y$ years is then obtained by multiplying the depreciation $\gamma$ of every
vehicle over that period with the VC required to cover all daily
trips, and then multiplying the distance cost with the average daily
VMT over the time period considered, assuming trips are only made on
weekdays.

\begin{figure}[!t]
	\centering
	\begin{minipage}{.49\textwidth}
		\centering
		\includegraphics[width=1.0\linewidth]{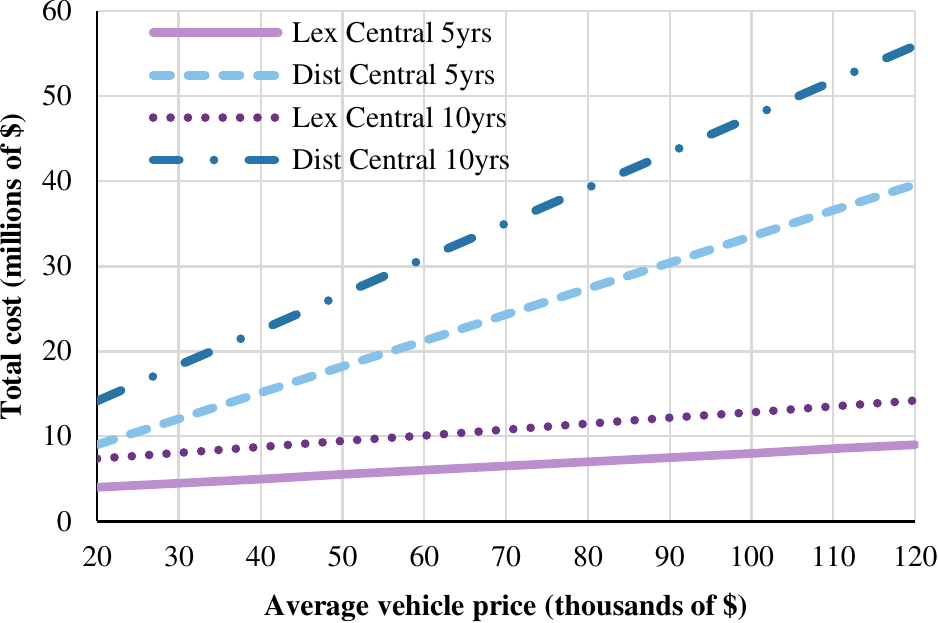}
		\caption{Total Vehicle Cost for CTSPAV Platform Inside City Limits Over 5 and 10 Years.}
		\label{fig:vehicle_cost_inside}
	\end{minipage}
	\hspace{0.0\textwidth}
	\begin{minipage}{.49\textwidth}
		\centering
		\includegraphics[width=1.0\linewidth]{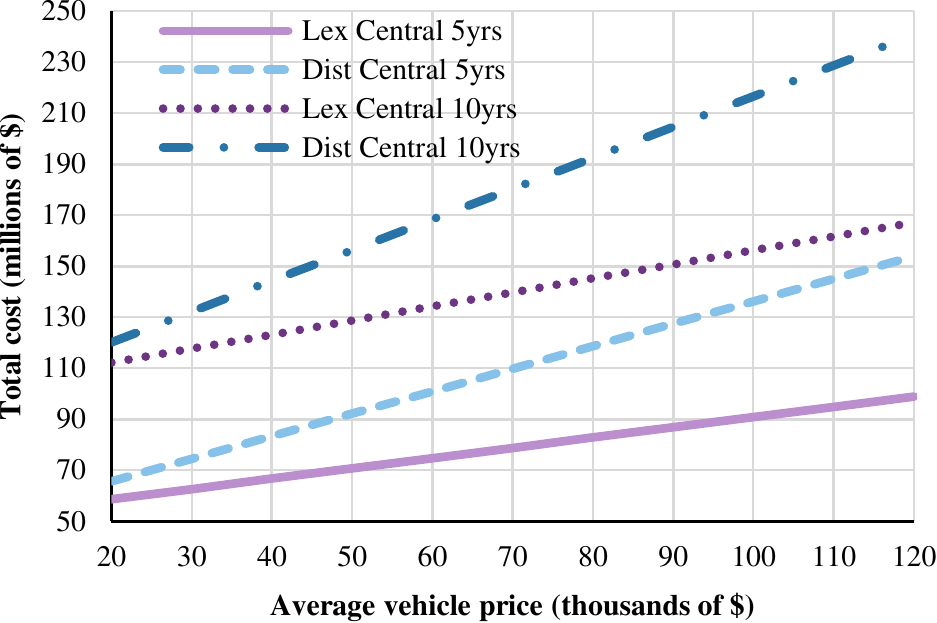}
		\caption{Total Vehicle Cost for CTSPAV Platform Outside City Limits Over 5 and 10 Years.}
		\label{fig:vehicle_cost_outside}
	\end{minipage}
\end{figure}

Figures \ref{fig:vehicle_cost_inside} and
\ref{fig:vehicle_cost_outside} show results of the vehicle cost
analysis over 5- and 10-years for clusters inside and outside city
limits respectively. In each figure, the horizontal axis displays the
range of possible initial prices of an AV, whereas the vertical axis
displays the corresponding total vehicle cost. The figures show that,
for the range of initial prices considered, the total vehicle cost is
always dominated by the contributions from the vehicle-age
cost. Therefore, the distance-minimization objective that requires
larger VCs is always more expensive than the lexicographic objective,
and this is true regardless of the location of the clusters or the
time period considered for the analysis.

The operating cost considers the average annual cost of a parking
permit for each vehicle combined with an estimated fixed cost for
installing a charging station for each AV. For this analysis, an
annual parking permit cost of \$800 is used together with a charging
station installation cost of \$1,400 per vehicle. The total
operating cost over $y$ years is then simply calculated by
multiplying the charging station installation cost with the VC
required to cover all daily trips, and then multiplying the annual
parking permit cost with the VC and the number of years.

\begin{figure}[!t]
	\centering
	\begin{minipage}{.49\textwidth}
		\centering
		\includegraphics[width=1.0\linewidth]{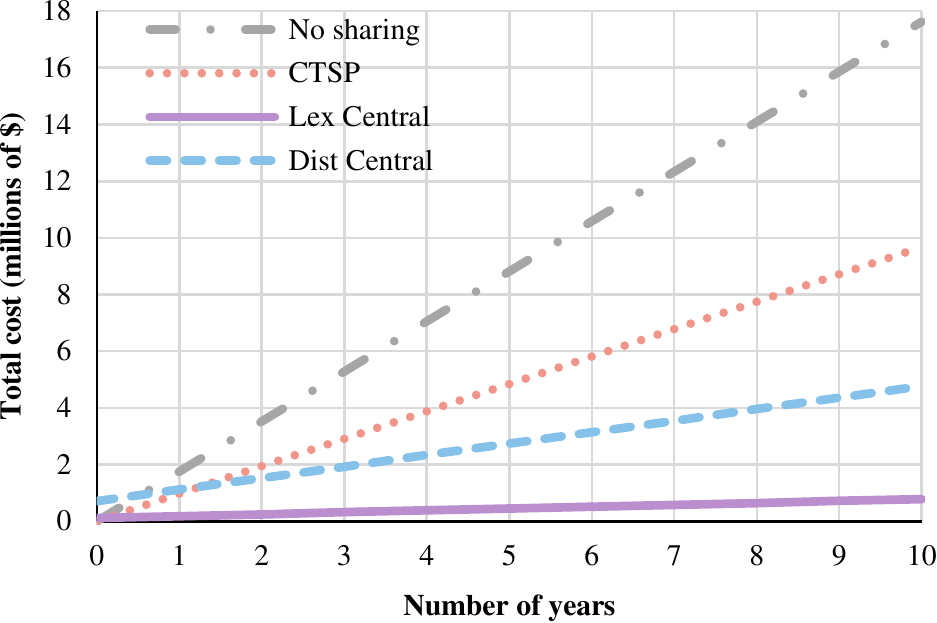}
		\caption{Total Operating Cost for CTSPAV Platform Inside City Limits Over 10 Years.}
		\label{fig:infra_cost_inside}
	\end{minipage}
	\hspace{0.0\textwidth}
	\begin{minipage}{.49\textwidth}
		\centering
		\includegraphics[width=1.0\linewidth]{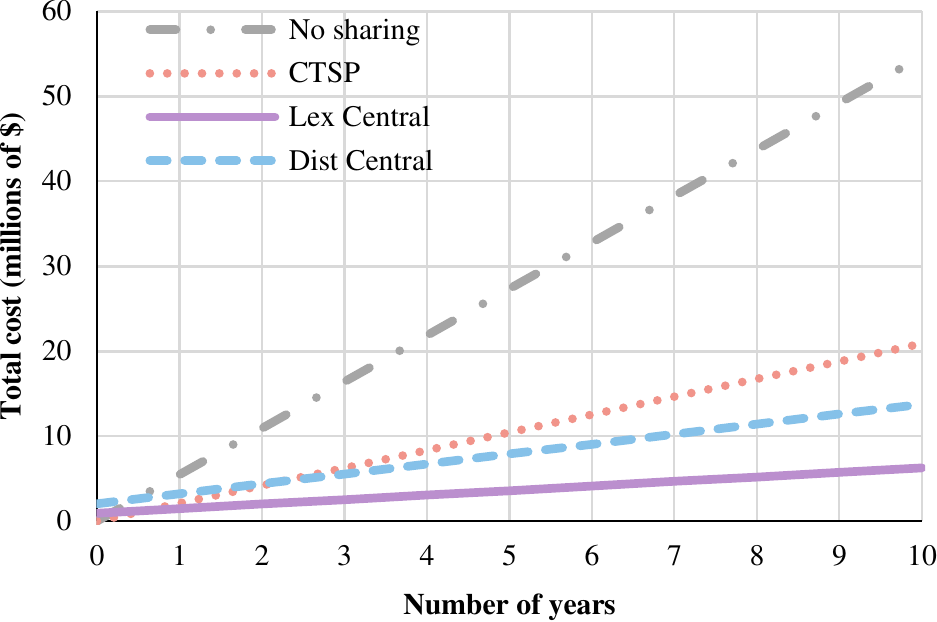}
		\caption{Total Operating Cost for CTSPAV Platform Outside City Limits Over 10 Years.}
		\label{fig:infra_cost_outside}
	\end{minipage}
\end{figure}

Figures \ref{fig:infra_cost_inside} and \ref{fig:infra_cost_outside}
display results of the total operating cost as a function of the
number of years for clusters inside and outside city limits
respectively. The total parking cost of the CTSP and of the vehicles
under no-sharing conditions are also included for additional
perspective. Since there are no fixed costs associated with the CTSP
or with the no-sharing conditions, their operating costs are
lower than those of the CTSPAV in the beginning. However, as time
increases, parking costs start to dominate, causing methods with
larger VCs to be more expensive. In fact, inside city limits, both
CTSPAV methods considered become cheaper as early as the second year,
while the same happens as early as the third year outside the
city. The gap in operating costs between the different methods
considered also increases with time and highlights the cost
effectiveness of the CTSPAV with the lexicographic objective as it
uses the fewest number of vehicles.

\subsection{Sensitivity to $\Delta$}

\begin{figure}[!t]
	\centering
	\includegraphics[width=1.0\linewidth]{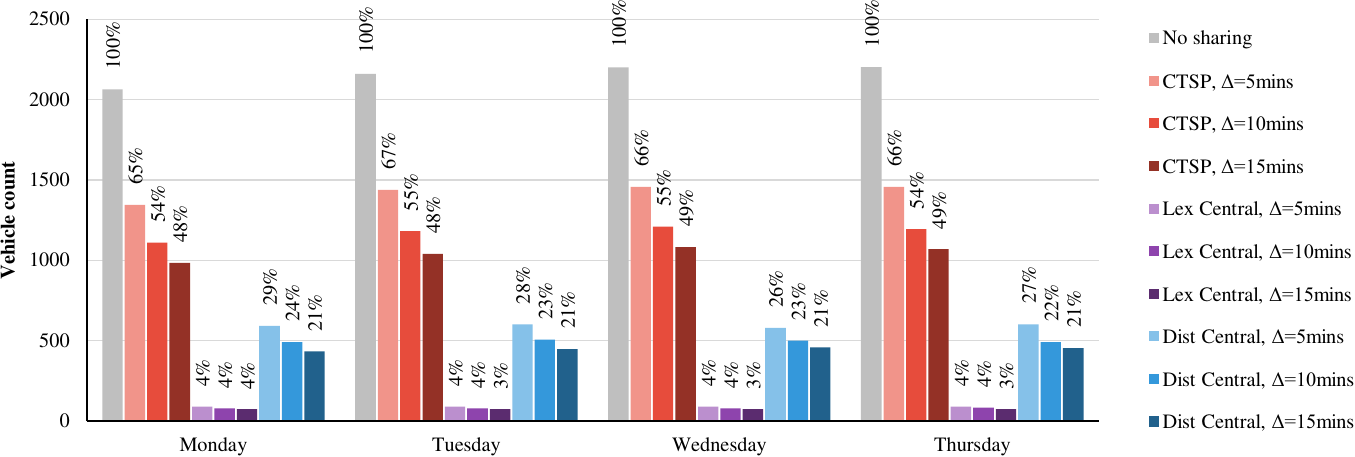}
	\caption{Aggregate Vehicle Count Results Inside City Limits for $\Delta = \{5, 10, 15\}$ mins.}
	\label{fig:vehcount_inside_deltasense}
\end{figure}

\begin{figure}[!t]
	\centering
	\includegraphics[width=1.0\linewidth]{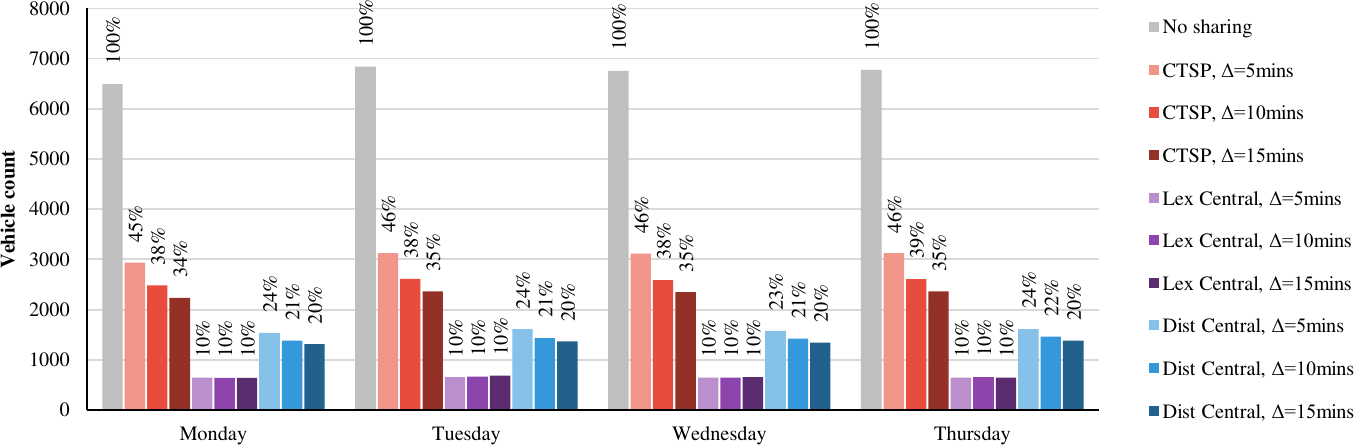}
	\caption{Aggregate Vehicle Count Results Outside City Limits for $\Delta = \{5, 10, 15\}$ mins.}
	\label{fig:vehcount_outside_deltasense}
\end{figure}

\begin{figure}[!t]
	\centering
	\includegraphics[width=1.0\linewidth]{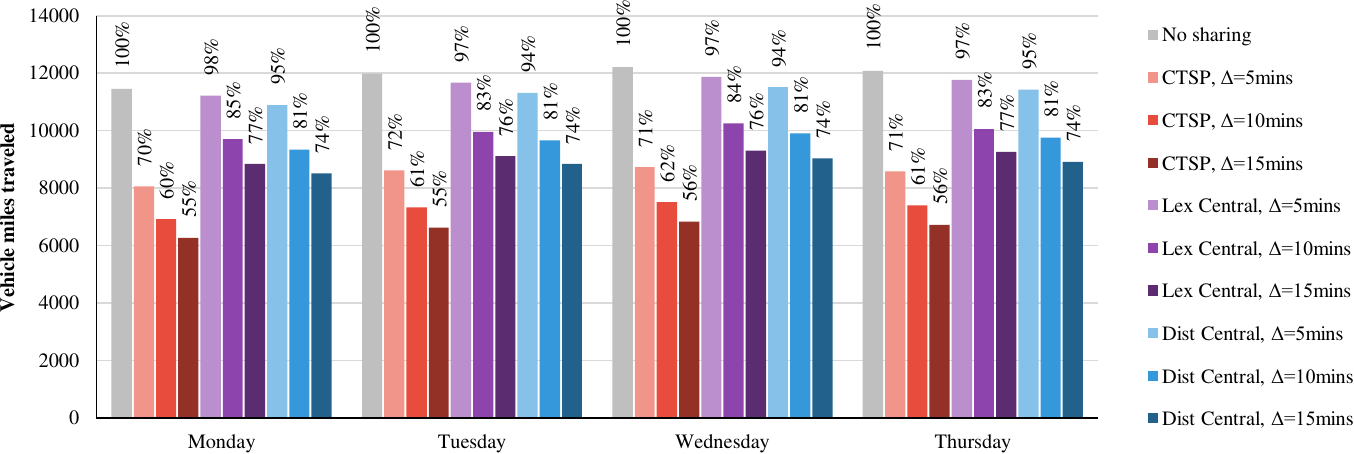}
	\caption{Aggregate Vehicle Miles Traveled Inside City Limits for $\Delta = \{5, 10, 15\}$ mins.}
	\label{fig:distance_inside_deltasense}
\end{figure}

\begin{figure}[!t]
	\centering
	\includegraphics[width=1.0\linewidth]{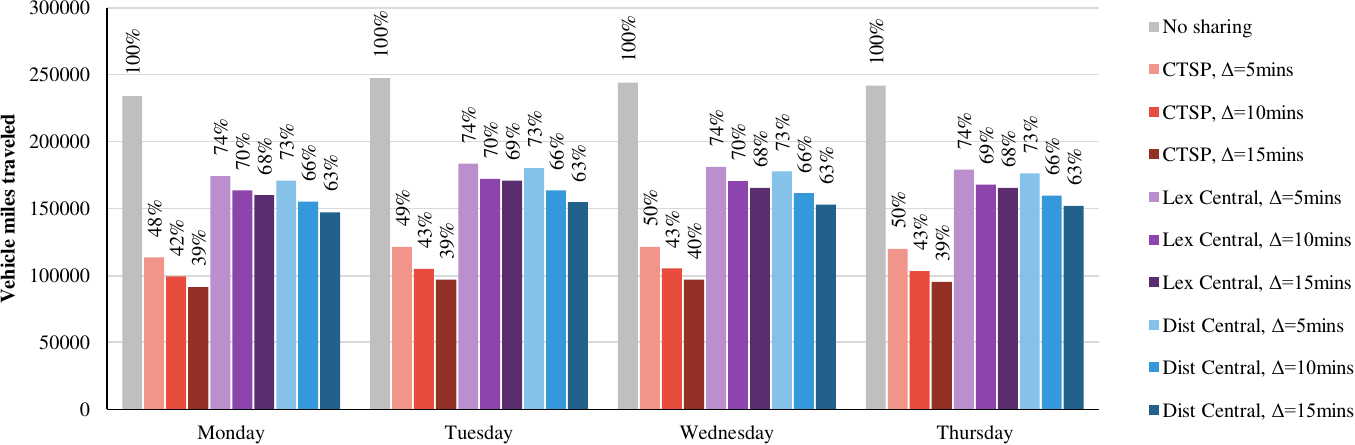}
	\caption{Aggregate Vehicle Miles Traveled Outside City Limits for $\Delta = \{5, 10, 15\}$ mins.}
	\label{fig:distance_outside_deltasense}
\end{figure}

The parameter $\Delta$ represents the maximum amount of time by which
each rider needs to shift (up or down) her desired arrival and
departure times at a parking structure. It has a direct
impact on the quality of service (QoS) of the riders, and it is
desirable to have $\Delta$ be as small as possible. Limiting its value
however restricts the flexibility of the schedules and could
negatively impact trip shareability. To study the impact of varying
$\Delta$ on the results of the CTSPAV, the procedure with the central
depot configuration is applied to optimize ridesharing for every
cluster with $\Delta$ set to $\{5, 15\}$ mins. The results are then
compared against those of $\Delta = 10$ mins. Figures
\ref{fig:vehcount_inside_deltasense} and
\ref{fig:vehcount_outside_deltasense} compare aggregated VCs of all
clusters inside and outside city limits respectively. Figures
\ref{fig:distance_inside_deltasense} and
\ref{fig:distance_outside_deltasense} then compare aggregated VMT from
all clusters inside and outside the city respectively. Sensitivity
results of the CTSP are also included in the figures, as well as the
results under no-sharing conditions and the percentage of each
quantity as a fraction of the no-sharing results for additional
perspective.

The results show that reducing $\Delta$ to 5 mins adversely affects
the vehicle reduction capability of the CTSP, reducing its VC by
approximately 12\% inside the city (resp. 8\% outside the
city). Increasing $\Delta$ to 15 mins improves vehicle reduction by
approximately 6\% inside the city (resp. 4\% outside the city). This
is the evidence of a trade-off between QoS and trip shareability. In
contrast, VC results of the CTSPAV with the lexicographic objective
exhibit very little sensitivity to $\Delta$, whereby the VCs change by
$\leq 1\%$ as $\Delta$ is varied by $\pm 5$ mins. \emph{This bodes
  very well for the CTSPAV as it indicates that a reduction in
  $\Delta$ to improve the quality of service for the riders will not
  have a significant impact on its vehicle reduction
  capability}. Finally, the VCs of the CTSPAV with the
distance-minimization objective display a modest sensitivity to
$\Delta$, whereby decreasing $\Delta$ by 5 mins degrades VCs by
approximately 5\% inside the city (resp. 3\% outside the city), and
increasing it by 5 mins improves VCs by approximately 2\% insde the
city (resp. 1\% outside the city).

The sensitivity analysis on aggregated VMT paints a different picture,
whereby every method considered exhibits comparable sensitivity to the
variations in $\Delta$. The results once again display a trade-off,
this time between QoS and travel distance reduction, evident from the
latter decreasing as $\Delta$ is increased and vice versa. The
increase in VMT of the CTSPAV, regardless of objective function or
position of clusters, when $\Delta$ is reduced stems from the
reduction of opportunities for trip aggregation as a result of the
tighter time windows. The AVs would therefore need to increase
back-and-forth traveling between the parking structures and
neighborhoods to serve the same amount of trips, leading to
corresponding increases in their travel distance.

\subsection{Sensitivity to $R$}

The parameter $R$ is yet another parameter that may influence the QoS
of the riders, as it directly influences the maximum amount of time
every rider spends in a vehicle. Limiting the value of $R$ improves QoS as
it leads to shorter ride durations. However, it also results in less
flexible trip schedules, which could consequently reduce the potential
for trip aggregation. Therefore, one would anticipate a trade-off
between QoS and trip shareability when varying $R$ similar to that
observed in the sensitivity analysis on $\Delta$. To investigate this
trade-off, a sensitivity analysis is conducted by setting $R = \{25\%,
75\%\}$, applying the procedure with the central depot configuration
to optimize ridesharing for every cluster, and comparing the results
with those of $R = 50\%$.

\begin{figure}[!t]
	\centering
	\includegraphics[width=1.0\linewidth]{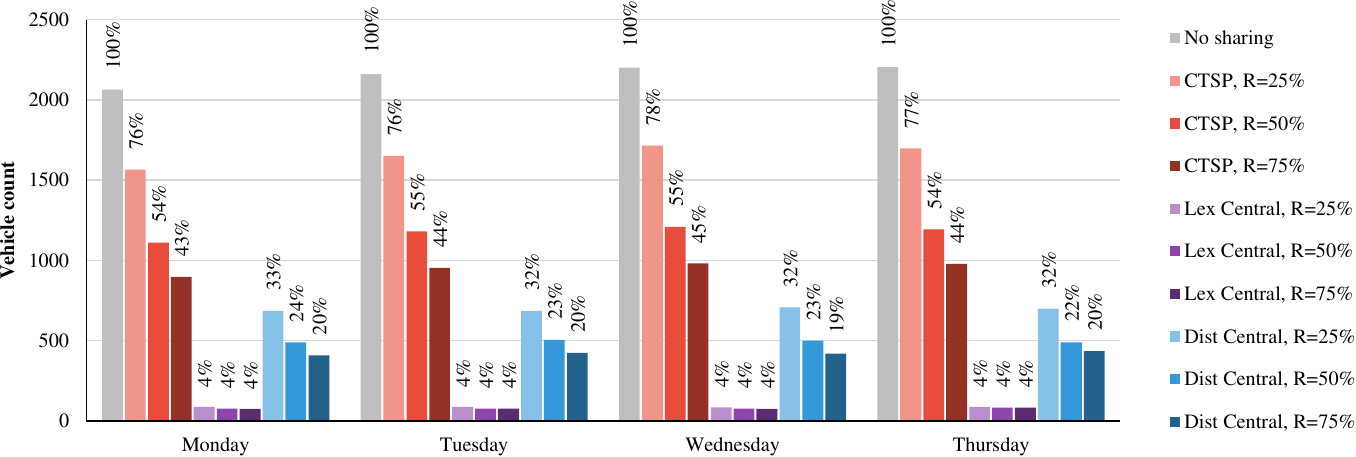}
	\caption{Aggregate Vehicle Count Results Inside City Limits for $R = \{25\%, 50\%, 75\%\}$.}
	\label{fig:vehcount_inside_Rsense}
\end{figure}

\begin{figure}[!t]
	\centering
	\includegraphics[width=1.0\linewidth]{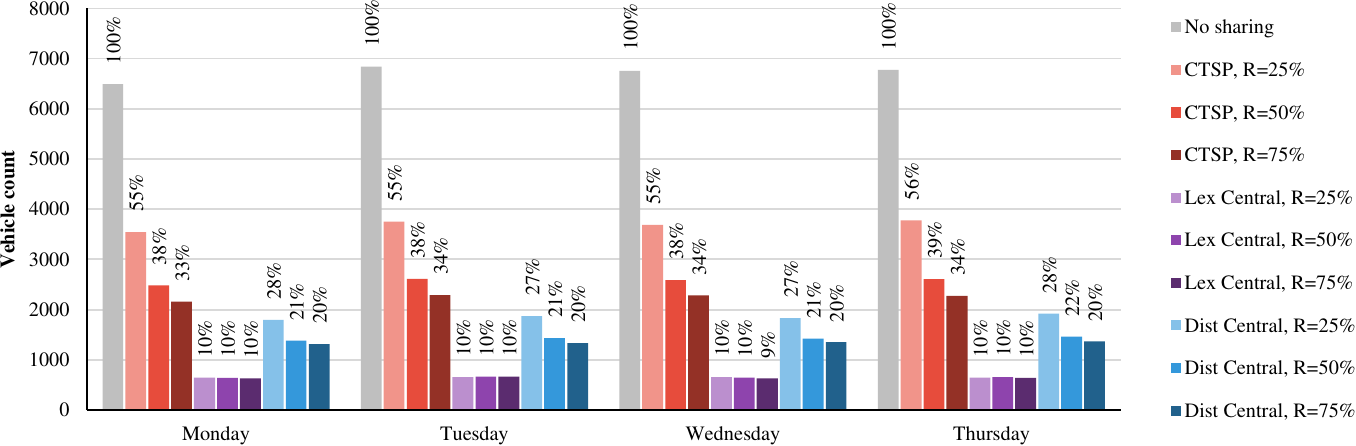}
	\caption{Aggregate Vehicle Count Results Outside City Limits for $R = \{25\%, 50\%, 75\%\}$.}
	\label{fig:vehcount_outside_Rsense}
\end{figure}

\begin{figure}[!t]
	\centering
	\includegraphics[width=1.0\linewidth]{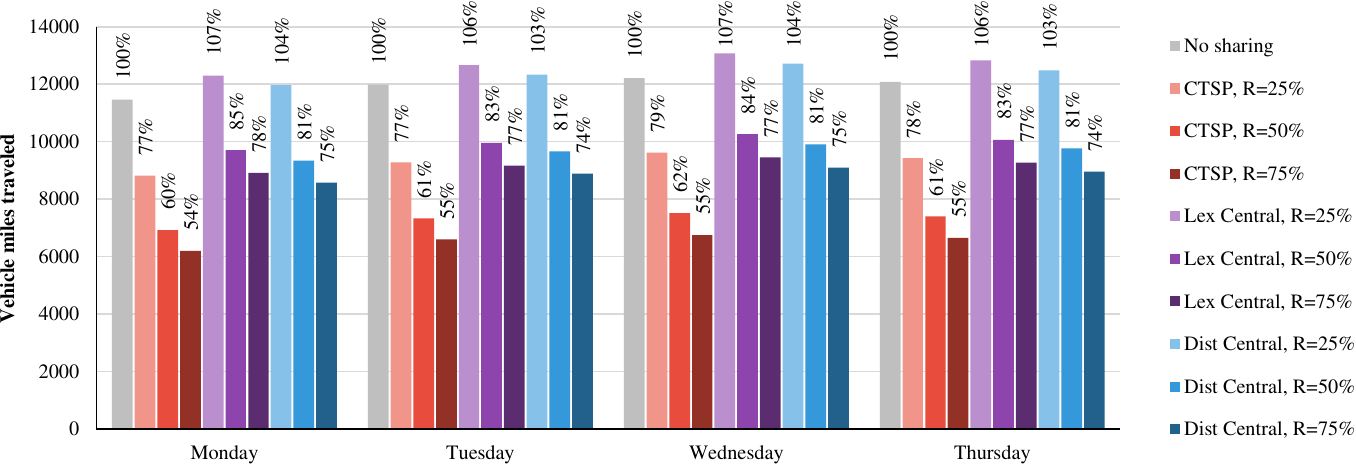}
	\caption{Aggregate Vehicle Miles Traveled Inside City Limits for $R = \{25\%, 50\%, 75\%\}$.}
	\label{fig:distance_inside_Rsense}
\end{figure}

\begin{figure}[!t]
	\centering
	\includegraphics[width=1.0\linewidth]{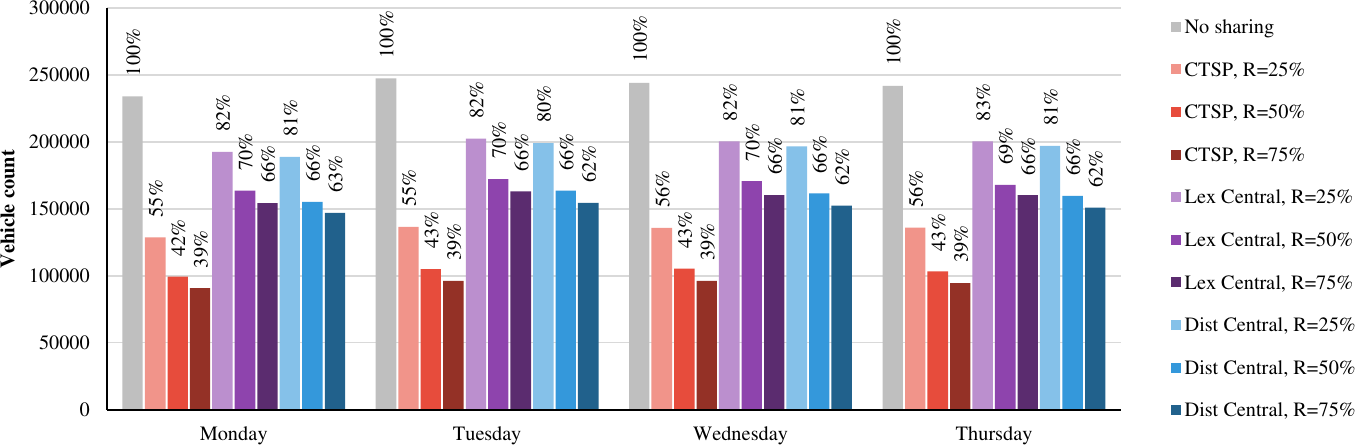}
	\caption{Aggregate Vehicle Miles Traveled Outside City Limits for $R = \{25\%, 50\%, 75\%\}$.}
	\label{fig:distance_outside_Rsense}
\end{figure}

Results of the analysis on VCs inside and outside city limits are
summarized in Figures \ref{fig:vehcount_inside_Rsense} and
\ref{fig:vehcount_outside_Rsense} respectively, whereas the analysis
on VMT inside and outside city limits are displayed in Figures
\ref{fig:distance_inside_Rsense} and \ref{fig:distance_outside_Rsense}
respectively. Similar to previous analyses, results from the original
CTSP, from the no-sharing condition, and percentages of each quantity
as a fraction of the no-sharing results are included for reference.

The VC results of the CTSP display the expected trade-off described
earlier where they appear to be very sensitive to variations in
$R$. Vehicle reduction significantly degrades by approximately 22\%
inside the city (resp. 17\% outside the city) when $R$ is reduced to
25\%, and it improves by approximately 11\% inside the city (resp. 5\%
outside the city) when $R$ is increased to 75\%. This indicates that
attempts to reduce the maximum ride durations will result in
significantly reduced trip shareability in the CTSP. In contrast, such
a drawback is not evident from the results of the CTSPAV with the
lexicographic objective, as its VCs exhibit very little sensitivity to
changes in $R$. The results show changes that are $<1\%$ both inside
and outside the city as $R$ is varied by $\pm
25\%$\footnotemark. \emph{This result provides another positive
  outlook for the CTSPAV as it promises that the vehicle reductions
  will not be adversely affected by attempts to increase QoS by
  reducing the maximum ride durations of the riders}. Finally, the VC
results of the CTSPAV with the distance minimization objective
displays a relatively modest sensitivity to $R$, whereby decreasing
$R$ to 25\% increased VCs by approximately 9\% inside the city
(resp. 6\% outside the city), and increasing $R$ to 75\% reduced VCs
by approximately 3\% inside the city (resp. 1\% outside the city).

\footnotetext{For the clusters outside city limits, the column-generation phase generated on average 34\% more columns for $R = 75\%$ than it did for $R = 50\%$, which caused the MIP to be significantly harder to solve. The MIP time limit for these instances is therefore extended to 2 hrs to account for this increase in complexity.}

Analysis of the VMT results reveals different observations. Inside city limits, both objective functions of the CTSPAV appear to be more sensitive to changes in $R$ than the CTSP; decreasing $R$ to 25\% leads to approximately a 23\% increase in VMT for both versions of the CTSPAV (compared to a 17\% increase for the CTSP). In fact, the increase is so significant in this case that their aggregated values exceed those for the no-sharing condition. This indicates that the opportunities for trip aggregation in this case is diminished to the point that any savings in travel distance is overshadowed by a more significant increase in back-and-forth traveling resulting from having to cover the trips with approximately the same amount of vehicles. Outside city limits, the VMT of every method considered exhibit comparable sensitivity to changes in $R$, once again displaying a trade-off between QoS and travel distance reduction.

\subsection{Effect of Increasing Vehicle Capacity}

\begin{figure}[!t]
	\centering
	\begin{minipage}{.49\textwidth}
		\centering
		\includegraphics[width=1.0\linewidth]{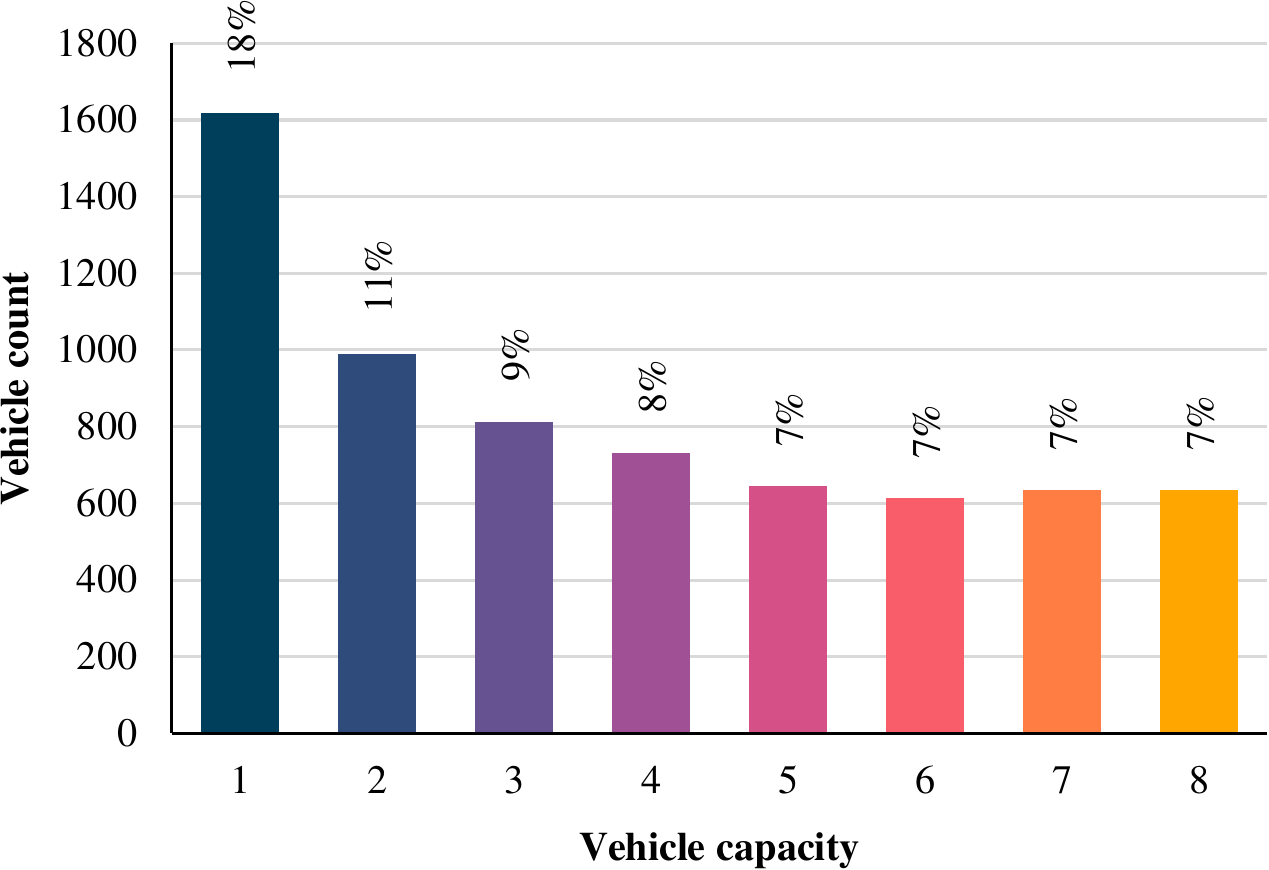}
		\caption{Effect of Increasing Vehicle Capacity on Aggregated Vehicle Count.}
		\label{fig:vc_k}
	\end{minipage}
	\hspace{0.0\textwidth}
	\begin{minipage}{.49\textwidth}
		\centering
		\includegraphics[width=1.0\linewidth]{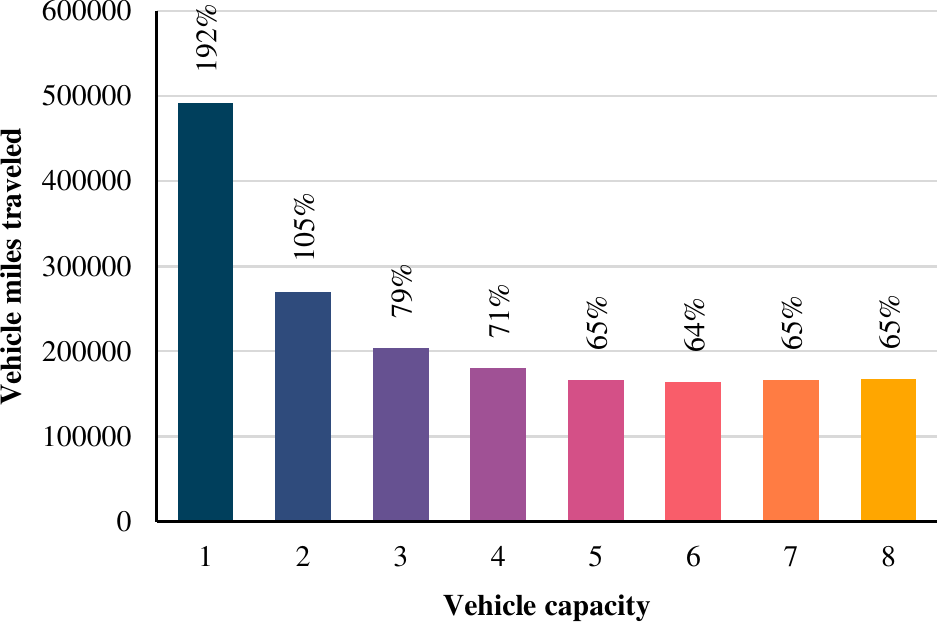}
		\caption{Effect of Increasing Vehicle Capacity on Aggregated Vehicle Miles Traveled.}
		\label{fig:vmt_k}
	\end{minipage}
\end{figure}

This section explores the effects of increasing vehicle capacity, $K$, on the results of the CTSPAV to investigate if there are any benefits to be gained from using vehicles with larger capacity. The results are obtained by varying $K$ between 1 to 8 and applying the CTSPAV procedure with the central depot configuration on all trips made on the Wednesday of week 2.\footnotemark Figure \ref{fig:vc_k} shows the aggregated vehicle counts from all clusters for every $K$ value, while Figure \ref{fig:vmt_k} does the same for the aggregated vehicle miles traveled. The percentages in both figures indicate each quantity as a fraction of the no-sharing values.

\footnotetext{To accommodate the exponential increase in problem complexity that results from the increase in $K$, the following changes were applied to the CTSPAV procedure for $K > 4$. A 10-minutes timeout is applied to the label-setting dynamic program of the PSP\textsubscript{CTSPAV}. The algorithm is terminated if it exceeds the timeout and it just returns all complete, non-dominated mini routes discovered at that point with negative reduced costs. The column-generation procedure is also seeded all with feasible mini routes covering up to 2 trips which are obtained from an exhaustive search procedure. Finally, the time budget for the column-generation phase is extended to 2 hours.}

Both figures reveal that marginal improvements (reductions) in VC and VMT decrease with increasing $K$. In fact, they show that almost no improvement is obtained for both VC and VMT beyond $K = 5$, and the biggest improvement is obtained when increasing $K$ from 1 to 2. These results are consistent with the findings from \cite{farhan2018} which focused on an SAV system for on-demand trips, whereby they discovered that the decreases to the fleet size was marginal for $K > 2$. It must also be noted that the marginal benefits obtained from increasing $K$ is accompanied by an exponential increase in problem complexity (in fact, the very slight increase in VC and VMT for $K > 6$ can be attributed the increase in problem complexity and the procedure not being able to find better solutions even with the extended time budget). The general trend of diminishing marginal benefits can probably be attributed to the unavailability of more trips with compatible itineraries that are needed to fully utilize the larger vehicle capacities from the problem instances considered.

\begin{figure}[!t]
	\centering
	\includegraphics[width=0.5\linewidth]{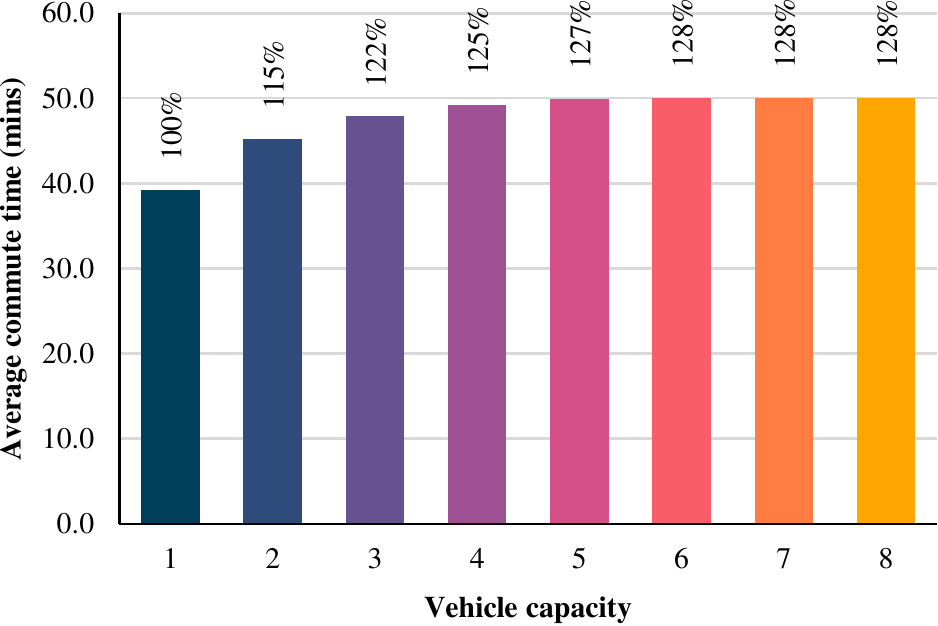}
	\caption{Effect of Increasing Vehicle Capacity on Average Commute Time.}
	\label{fig:act_k}
\end{figure}

Finally, Figure \ref{fig:act_k} shows the effect of increasing $K$ on the average ride duration per commuter. The percentages represent each quantity as a fraction of the value when $K = 1$. Unsurprisingly, the best average commute time is obtained when $K = 1$ as the rides are not shared under this setting, and the times increase with increasing $K$ due to the corresponding increase in ridesharing. The diminishing nature of the marginal increases in average commute time as $K$ is increased also appears to mirror the marginal decreases in VC or VMT. The figure also shows that when $K = 4$, the average commute time only increases by 25\% (even though $R = 0.50$), which bodes very well for maintaining the QoS of riders of the CTSPAV.

\section{Conclusion}
\label{sec:conclusion}

This work originated from the desire to understand the potential
benefits of autonomous vehicles on a car-pooling platform that
maximizes ridesharing for the commuting trips of a community. The
central problem powering the platform is the CTSP: it was 
originally conceived by \citet{hasan2018} to
reduce parking utilization and traffic congestion in urban areas by
leveraging the structure of commuting patterns and urban
communities. The CTSP was shown to reduce the number of vehicles
significantly on a real case study. 

Given that vehicles are idle for most of the day by definition of
car-pooling, it is interesting to study how much additional reduction
in fleet size would come from using autonomous vehicles, as well as
the impact on miles traveled. To answer that question, this paper
defined the CTSPAV problem and proposed two column-generation
procedures to obtain high-quality solutions. The first approach (the
CTSPAV procedure) assembles feasible mini routes into an overall
routing plan, while the second approach (the DARP procedure) reduces
the CTSPAV into a DARP. The optimization problems considered (1) a
lexicographic objective that first minimizes the required vehicle
count and then their total distance, and (2) a distance minimization
objective that just minimizes the total travel distance. 

The CTSPAV was evaluated on a large-scale, real-world dataset of
commute trips from the city of Ann Arbor, Michigan, which contains
detailed information for an average of 9,000 daily trips over a month.
The experimental results revealed that the CTSPAV procedure with the
lexicographic objective reduces daily vehicle usage by 92\% while at
the same time also reducing vehicle miles traveled by 30\%. Its
vehicle reduction results represent a 34\% improvement relative to
that of the CTSP. Examining the solutions show that the CTSPAV
generates significantly longer routes by traveling back and forth
between the communities and the commuting destination (and vice
versa). A cost analysis also showed that fleets of AV vehicles are
eminently viable from an economic standpoint in this setting. Finally,
sensitivity analyses revealed that the vehicle count results of the
lexicographic objective are more resilient to variations in the size
of the time windows and the length of the ride-duration limits of the
trips. The number of vehicles required to cover all trips changed by
less than 1\% as the two parameters are varied, while the same cannot
be said about the CTSP.

\section*{Acknowledgement}

We would like to thank Stephen Dolen from Logistics, Transportation,
and Parking of the University of Michigan for his assistance in
obtaining the dataset used in this research. Part of this research was
funded by the Rackham Graduate Student Research Grant, computational resources and services provided by Advanced Research Computing at the University of
Michigan, and NSF Leap HI proposal NSF-1854684.





\bibliographystyle{elsarticle-harv} 
\bibliography{references}

\newpage
\appendix

\section{The Pricing Subproblem PSP\textsubscript{CTSPAV}}
\label{appendix:pricing}

\paragraph{Construction of graphs for PSP\textsubscript{CTSPAV}}

Let $v^i_t$ denote a virtual sink node for graph $\mathcal{G}_i^+$ or
$\mathcal{G}_i^-$. Since a mini route covers only inbound trips or
outbound trips, $\mathcal{G}_i^+= (\mathcal{N}_i^+,\mathcal{A}_i^+)$
(resp. $\mathcal{G}_i^-= (\mathcal{N}_i^-,\mathcal{A}_i^-)$) contains,
in addition to $v^i_t$, only the nodes for inbound trips
(resp. outbound trips), i.e., $\mathcal{N}_i^+ =
\mathcal{P}^+\cup\mathcal{D}^+\cup\{v^i_t\}$ (resp. $\mathcal{N}_i^- =
\mathcal{P}^-\cup\mathcal{D}^-\cup\{v^i_t\}$). The set
$\mathcal{A}_i^+$ (resp. $\mathcal{A}_i^-$) then represents all
feasible edges for $\mathcal{G}_i^+$ (resp. $\mathcal{G}_i^-$), i.e.,
location pairs from $\mathcal{N}_i^+$ (resp. $\mathcal{N}_i^-$) that
satisfy a priori route feasibility constraints. Without loss of
generality, the following elaborates further on how $\mathcal{G}_i^+$
is constructed.

Construction of $\mathcal{G}_i^+=(\mathcal{N}_i^+,\mathcal{A}_i^+)$
begins with the introduction of the set of nodes $\mathcal{N}_i^+ =
\mathcal{P}^+\cup\mathcal{D}^+\cup\{v^i_t\}$ and a set of
fully-connected edges
$\mathcal{A}_i^+=\{(u,v):u,v\in\mathcal{N}_i^+,u\neq v\}$. A
ride-duration limit $L_u$ is associated with each node
$u\in\mathcal{P}^+$, a time window $[a_u,b_u]$ and service duration
$s_u$ are associated with each node
$u\in\mathcal{P}^+\cup\mathcal{D}^+$, and a travel time $\tau_{(u,v)}$
and a reduced cost $\bar{c}_{(u,v)}$ are associated with each edge
$(u,v)\in\mathcal{A}^+_i$. As the goal is to find a feasible mini
route from $i$ to $v^i_t$ with minimum reduced cost, $\bar{c}_{(u,v)}$
is defined as follows so that the total cost of any path from $i$ to
$v^i_t$ is equivalent to that defined in \eqref{mini_route_rc}.
\begin{equation}\label{edge_rc}
\bar{c}_{(u,v)} = 
\begin{cases}
-\pi_u - \mu_{(u,v)}&\qquad \forall (u,v)\in\mathcal{A}^+_i\setminus\delta^-(v^i_t):u\in\mathcal{P}^+\\
-\mu_{(u,v)}&\qquad \forall (u,v)\in\mathcal{A}^+_i\setminus\delta^-(v^i_t):u\in\mathcal{D}^+\\
0&\qquad \forall (u,v)\in\delta^-(v^i_t)
\end{cases}
\end{equation}

\noindent
Edges from $\mathcal{A}^+_i$ that cannot belong to any feasible mini
route are then identified by pre-processing time-window, pairing,
precedence, and ride-duration limit constraints. Prior to this
pre-processing step, knowledge of $i$ being the source of the path
sought from $\mathcal{G}_i^+$ allows the time windows of all nodes
$u\in\mathcal{P}^+\cup\mathcal{D}^+$ to be tightened, by sequentially
increasing their lower bounds using the following rules proposed by
\citet{dumas1991}:
\begin{itemize}
\item $a_u = \max \{a_u, a_i+s_i+\tau_{(i,u)}\}, \forall u\in\mathcal{P}^+\setminus\{i\}$	
\item $a_{n+u} = \max \{a_{n+u}, a_u+s_u+\tau_{(u,n+u)}\}, \forall u\in\mathcal{P}^+\setminus\{i\}$	
\end{itemize}

\begin{figure}[!t]
	\centering
	\includegraphics[width=0.5\linewidth]{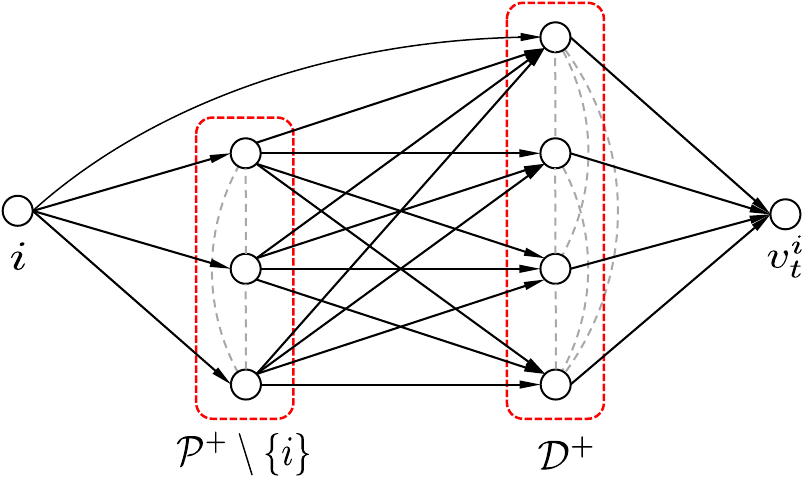}
	\caption{Graph $\mathcal{G}^+_i$ (Each Dotted Line Represents a Pair of Bidirectional Edges).}
	\label{fig:graph2}
\end{figure}

\noindent
The following sets of infeasible edges are then identified and consequently removed from $\mathcal{A}^+_i$:
\begin{enumerate}[(a)]
\item Direct trips to source $i$ and from sink $v^i_t$: $\delta^-(i)\cup\delta^+(v^i_t)$
\item Precedence of pickup and drop-off nodes: 	
\begin{itemize}
	\item $\{(i,v) : v\in\mathcal{D}^+\setminus\{n+i\}\}$
	\item $\{(u,v^i_t) : u\in\mathcal{P}^+\}$
	\item $\{(u,v) : u\in\mathcal{D}^+\wedge v\in\mathcal{P}^+\}$
\end{itemize}
\item Time windows along each edge: $\{(u,v):(u,v)\in\mathcal{A}^+_i \setminus \delta^-(v^i_t) \wedge a_u + s_u + \tau_{(u,v)} > b_v\}$
\item Ride-duration limit of each commuter: $\{(u,v),(v,n+u):u\in\mathcal{P}^+ \wedge v\in\mathcal{P}^+\cup\mathcal{D}^+ \wedge u\neq v \wedge \tau_{(u,v)} + s_v + \tau_{(v,n+u)} > L_u\}$
\item Time windows and ride-duration limits of pairs of trips:
\begin{itemize}
	\item $\{(u,n+v):u,v\in\mathcal{P}^+ \wedge u\neq v \wedge \neg feasible(v\rightarrow u\rightarrow n+v\rightarrow n+u)\}$
	\item $\{(n+u,v):u,v\in\mathcal{P}^+ \wedge u\neq v \wedge \neg feasible(u\rightarrow n+u\rightarrow v\rightarrow n+v)\}$
	\item $\{(u,v):u,v\in\mathcal{P}^+ \wedge u\neq v \wedge \neg feasible(u\rightarrow v\rightarrow n+u\rightarrow n+v) \wedge \neg feasible(u\rightarrow v\rightarrow n+v\rightarrow n+u)\}$
	\item $\{(n+u,n+v):u,v\in\mathcal{P}^+ \wedge u\neq v \wedge \neg feasible(u\rightarrow v\rightarrow n+u\rightarrow n+v) \wedge \neg feasible(v\rightarrow u\rightarrow n+u\rightarrow n+v)\}$
\end{itemize}
\end{enumerate}
An example of graph $\mathcal{G}_i^+$ that results from the removal of the infeasible edges is shown in Figure \ref{fig:graph2}.

\paragraph{Search for path with minimum reduced cost}

The final step in PSP\textsubscript{CTSPAV} is to find a path from $i$
to $v^i_t$ from each graph $\mathcal{G}_i^+:i\in\mathcal{P}^+$ and
$\mathcal{G}_i^-:i\in\mathcal{P}^-$. On top of having the minimum
cost, the path must also represent a \emph{feasible} mini
route. Notice that, as shown in Figure \ref{fig:graph2}, construction
of the graphs eliminates all edges $\{(u,v) : u\in\mathcal{D}\wedge
v\in\mathcal{P}\}$. Therefore, any path from $i$ to $v^i_t$ is
guaranteed to begin with a pickup phase that only visits nodes in
$\mathcal{P}$ and end with a drop-off phase that only visits nodes in
$\mathcal{D}$. In other words, the graphs ensure the precedence
feasibility constraint, which requires all pickup nodes of a mini
route to precede all of its drop-off nodes, is satisfied by
construction. All that remains to guarantee the feasibility of any
path found is to ensure that the path:
\begin{enumerate}[1.]
	\item Visits each node within its specified time window,
	\item Visits the corresponding drop-off nodes of every pickup node visited,
	\item Satisfies the ride-duration limit of each rider served,
	\item Respects the vehicle capacity, and
	\item Visits each node at most once (i.e., it is elementary).
\end{enumerate}
  
The problem of finding a feasible, least-cost path from $i$ to $v^i_t$
is essentially an Elementary Shortest Path Problem with Resource
Constraints (ESPPRC) \citep{irnich2005}, whereby the remaining
route-feasibility constraints can be modeled with resource
constraints. In fact, this ESPPRC is identical to that in the
column-generation pricing problem of the DARP considered in
\citet{gschwind2015}. They proposed a dynamic-programming,
label-setting algorithm that uses strong label dominance rules to find
the optimal solution to the problem.

While their algorithm can also be used to effectively solve this
work's ESPPRC, the algorithm has a pre-requisite; it requires that the
edge reduced costs $\bar{c}_{(u,v)}$ satisfy the \emph{delivery
  triangle inequality} (DTI). The DTI, introduced by
\citet{ropke2009}, requires that the reduced costs fulfill
$\bar{c}_{(u,v)} \leq \bar{c}_{(u,w)} + \bar{c}_{(w,v)}$ for all edges
$(u,v),(u,w),(w,v)\in\{\mathcal{A}^+_i:i\in\mathcal{P}^+\}\cup\{\mathcal{A}^-_i:i\in\mathcal{P}^-\}$
and $w\in\mathcal{D}$. Unfortunately, the reduced costs in this
problem do not satisfy the DTI. A cost-matrix transformation, also
proposed by \citet{ropke2009}, is therefore first applied to the
reduced costs to transform them into an equivalent set of costs that
does satisfy the DTI, after which the label-setting algorithm of
\citet{gschwind2015} is applied on each graph
$\mathcal{G}_i^+:i\in\mathcal{P}^+$ and
$\mathcal{G}_i^-:i\in\mathcal{P}^-$ to find the feasible, least-cost
paths from $i$ to $v^i_t$. Finally, note that enforcement of pairing
and precedence resource constraints of the label-setting algorithm on
the origin and destination node pairs of the graphs are sufficient to
ensure elementarity of the paths produced, as the graphs lack edges
$\{(u,v) : u\in\mathcal{D}\wedge v\in\mathcal{P}\}$ which are
necessary to produce cycles in the presence of the pairing and
precedence constraints. Therefore, additional resource constraints
dedicated specifically to ensure elementarity of the paths are not
necessary for PSP\textsubscript{CTSPAV}.

\section{Optimality Gaps and Computation Times}
\label{sec:gaps}

The optimality gaps and computation times of all experiments for the
CTSPAV procedure are summarized in Tables \ref{tab:gap_time_lex} and
\ref{tab:gap_time_dist}. In Table \ref{tab:gap_time_lex}, which
summarizes the results of problem instances for the lexicographic
objective, the first two columns specify the location of the clusters
and the configuration of their depots. The following two columns list
average values of the absolute gap for the vehicle count. The absolute
gap is calculated by taking the vehicle count results of the MIP and
subtracting from it its lower bound. The first gap uses results the
from RMP\textsubscript{CTSPAV}. Letting $Y_e^*$ be the value of $Y_e$
from RMP\textsubscript{CTSPAV} at convergence, the primal lower bound
to the vehicle count is given by $\ceil{\sum_{e\in\delta(v_s)}
  Y_e^*}$. The second gap uses the primary objective results of
RMP\textsubscript{DARP}, particularly $\ceil{z_\text{LB}'}$ as the
vehicle count lower bound. It was found to consistently provide a
stronger lower bound than the former, and therefore it is included
here to provide an additional perspective. The next column shows the
average optimality gap, which is given by $(z_\text{MIP} - z^*) /
z_\text{MIP}$ for each instance, where $z_\text{MIP}$ denotes the
MIP's final objective value, and $z^*$, the objective value of
RMP\textsubscript{CTSPAV} at convergence, provides a primal lower
bound for $z_\text{MIP}$. However, for a few problem instances
representing clusters outside city limits, the column-generation phase
did not converge within the time limit. For these instances, the dual
lower bound $z_\text{LB}$ is used in place of $z^*$ when calculating
the optimality gap, and the vehicle count lower bound from the CTSPAV
procedure is obtained by considering only the fixed cost contributions
to $z_\text{LB}$. All uncertainties are represented by the standard
error of the mean.

While the average optimality gap for the lexicographic objective is relatively high at approximately 70\% (resp. 76\%) for clusters inside the city (resp. outside the city), the vehicle count gap, which depicts the absolute gap of the primary objective, paints a different picture, averaging at 2.5 and 7.0 vehicles for clusters inside and outside city limits respectively. Moreover, when the lower bounds produced by the DARP procedure are used, the vehicle count gap averages at even lower values of 1.1 and 4.2 vehicles for clusters inside and outside city limits respectively. These lower values highlight a key strength of the DARP procedure, \emph{it consistently produces stronger lower bounds for the primary objective}. When these stronger lower bounds are used, some problem instances produced an absolute vehicle count gap of zero, indicating that their primary objective results are optimal.

\begin{table}[!t]
	\caption{Optimality Gaps and Computation Times of the CTSPAV Procedure with the Lexicographic Objective}
	\label{tab:gap_time_lex}
	\resizebox{\textwidth}{!}{
	\begin{tabular}{ccccccccc}
		\noalign{\rule{0pt}{10pt}}\hline\noalign{\smallskip}
		\begin{tabular}[c]{@{}c@{}}Cluster\\location\end{tabular} &
		\begin{tabular}[c]{@{}c@{}}Depot\\config.\end{tabular} &
		\multicolumn{1}{c}{\begin{tabular}[c]{@{}c@{}}Average\\vehicle\\count gap\end{tabular}} &
		\multicolumn{1}{c}{\begin{tabular}[c]{@{}c@{}}Average\\vehicle\\count gap\\(DARP LB)\end{tabular}} &
		\multicolumn{1}{c}{\begin{tabular}[c]{@{}c@{}}Average\\optimality\\gap (\%)\end{tabular}} &
		\multicolumn{1}{c}{\begin{tabular}[c]{@{}c@{}}Average\\mini route\\count\end{tabular}} &
		\multicolumn{1}{c}{\begin{tabular}[c]{@{}c@{}}Average\\colgen\\time (s)\end{tabular}} &
		\multicolumn{1}{c}{\begin{tabular}[c]{@{}c@{}}\% $\geq$\\colgen\\time\\limit\end{tabular}} &
		\multicolumn{1}{c}{\begin{tabular}[c]{@{}c@{}}\% $\geq$\\MIP\\time\\limit\end{tabular}} \\
		\noalign{\smallskip}\hline\noalign{\smallskip}
		\multirow{2}{*}{\begin{tabular}[c]{@{}c@{}}Inside\\ city\end{tabular}}  & Central & 2.55 $\pm$ 0.14 & 1.09 $\pm$ 0.15 & 70.5 $\pm$ 1.4 & 5394 $\pm$ 741 & 86 $\pm$ 28   & 0.0  & 100.0 \\
		\noalign{\smallskip}\cline{2-9}\noalign{\smallskip}
		& Local   & 2.50 $\pm$ 0.16 & 1.09 $\pm$ 0.16 & 69.8 $\pm$ 1.7 & 5699 $\pm$ 893 & 161 $\pm$ 84  & 0.0  & 100.0 \\
		\noalign{\smallskip}\hline\noalign{\smallskip}
		\multirow{2}{*}{\begin{tabular}[c]{@{}c@{}}Outside\\ city\end{tabular}} & Central & 7.07 $\pm$ 0.45 & 4.19 $\pm$ 0.28 & 75.5 $\pm$ 2.1 & 19949 $\pm$ 1021 & 926 $\pm$ 141 & 32.4 & 100.0 \\
		\noalign{\smallskip}\cline{2-9}\noalign{\smallskip}
		& Local   & 7.03 $\pm$ 0.47 & 4.12 $\pm$ 0.30 & 75.5 $\pm$ 2.2 & 19631 $\pm$ 974 & 922 $\pm$ 133 & 36.8 & 100.0 \\
		\noalign{\smallskip}\hline\noalign{\smallskip}
	\end{tabular}
}
\end{table}

The next column shows the average number of mini routes, i.e. columns, generated in the column-generation phase, while the following two columns summarize the time spent in this phase. The latter of the two summarizes the fraction of problem instances in which the column-generation phase met or exceeded the 1 hour time limit (i.e., the fraction of instances whereby the column-generation phase did not converge within the time limit), while the former shows the average wall time spent on column generation \emph{for problem instances that did not exceed the time limit}. It can be seen that on average, the time spent on column generation is significantly higher on clusters outside city limits. On top of that, more than 30\% of the instances did not achieve convergence within the time limit. This can be attributed to the significantly larger number of columns generated from these clusters. Finally, the last column shows the fraction of problem instances in which the MIP solving phase met or exceeded its 1 hour time budget. For the lexicographic objective, the time limit was exceeded on all problem instances.

\begin{table}[!t]
	\caption{Optimality Gaps and Computation Times of the CTSPAV Procedure with the Distance Minimization Objective}
	\label{tab:gap_time_dist}
	\resizebox{\textwidth}{!}{
	\begin{tabular}{ccccccc}
		\noalign{\rule{0pt}{10pt}}\hline\noalign{\smallskip}
		\begin{tabular}[c]{@{}c@{}}Cluster\\ location\end{tabular} &
		\begin{tabular}[c]{@{}c@{}}Depot\\ config.\end{tabular} &
		\multicolumn{1}{c}{\begin{tabular}[c]{@{}c@{}}Average\\ optimality\\ gap (\%)\end{tabular}} &
		\multicolumn{1}{c}{\begin{tabular}[c]{@{}c@{}}Average\\ mini route\\ count\end{tabular}} &
		\multicolumn{1}{c}{\begin{tabular}[c]{@{}c@{}}Average\\ colgen\\ time (s)\end{tabular}} &
		\multicolumn{1}{c}{\begin{tabular}[c]{@{}c@{}}\% $\geq$ colgen \\ time limit\end{tabular}} &
		\multicolumn{1}{c}{\begin{tabular}[c]{@{}c@{}}\% $\geq$ MIP \\ time limit\end{tabular}} \\
		\noalign{\smallskip}\hline\noalign{\smallskip}
		\multirow{2}{*}{\begin{tabular}[c]{@{}c@{}}Inside\\ city\end{tabular}}  & Central & 1.75 $\pm$ 0.17 & 5396 $\pm$ 752   & 73 $\pm$ 22    & 0.0  & 36.4 \\
		\noalign{\smallskip}\cline{2-7}\noalign{\smallskip}
		& Local   & 1.35 $\pm$ 0.17 & 5453 $\pm$ 803   & 105 $\pm$ 37   & 0.0  & 13.6 \\
		\noalign{\smallskip}\hline\noalign{\smallskip}
		\multirow{2}{*}{\begin{tabular}[c]{@{}c@{}}Outside\\ city\end{tabular}} & Central & 3.05 $\pm$ 0.48 & 20533 $\pm$ 1147 & 960 $\pm$ 112  & 8.8  & 82.4 \\
		\noalign{\smallskip}\cline{2-7}\noalign{\smallskip}
		& Local   & 4.75 $\pm$ 1.04 & 21114 $\pm$ 1140 & 1126 $\pm$ 129 & 19.1 & 73.5 \\
		\noalign{\smallskip}\hline\noalign{\smallskip}
	\end{tabular}
}
\end{table}

Table \ref{tab:gap_time_dist} summarizes the results of problem instances for the distance-minimization objective. It displays the same information as Table \ref{tab:gap_time_lex}, except that it does not list the average vehicle count gaps. This is due to the lower bound of the vehicle count not being available from the objective. It can be seen that the optimality gap for the distance-minimization objective is excellent, being less than 2\% on average for cases inside the city and less than 5\% on average for cases outside. Similar to instances for the lexicographic objective, the average number of columns generated from clusters outside city limits is significantly larger than that from clusters inside, which consequently results in the significantly longer average time spent on the column-generation phase. Finally, the percentage of problem instances that exceed the time limits of the column-generation and MIP phases is fewer than that for the lexicographic objective in every case, and this can be attributed to the relatively stronger primal lower bound provided by the LP relaxation of MP\textsubscript{CTSPAV} when the distance-minimization objective is utilized. 

Table \ref{tab:darp_gap_time_lex} summarizes the average optimality gaps and computations times of the same set of problem instances for the DARP procedure with the lexicographic objective. Once again, all uncertainties are represented by the standard error of the mean. Its first two columns specify the cluster location and depot configuration of the instances. The next lists the average number of columns generated, followed by two columns which show the average absolute gap and average optimality gap of \emph{the primary objective value}. Letting $z_\text{MIP}'$ denote the final primary objective value of the MIP and recalling that $z_\text{LB}'$ denotes its dual lower bound, the absolute gap of an instance given by $z_\text{MIP}' - \ceil{z_\text{LB}'}$ while its optimality gap is given by $(z_\text{MIP}' - \ceil{z_\text{LB}'})/z_\text{MIP}'$. The average optimality gaps are relatively high at approximately 48\%. However, the absolute gaps reveal a different story, whereby they average at only 2.6 and 4.8 vehicles inside and outside city limits respectively. These values, however, are higher than those of the CTSPAV procedure which utilize the same lower bounds ($\ceil{z_\text{LB}'}$). 

The next column specifies the fraction of problem instances whereby the column-generation phase for the primary objective did not converge within the time limit, followed by one which specifies the fraction of instances whereby the column-generation phase for both objectives did not converge. When only the primary objective is considered, column-generation for more than two-thirds of the instances inside the city and almost all instances outside the city did not converge. On top of that, out of the 180 problem instances in total, only two had their column-generation phases for both objectives converge within the allocated time limit, even after more time has been allocated for this phase (recall 1.5 hours is allocated for this phase of the DARP procedure as opposed to only 1 hour for the CTSPAV procedure). This highlights the harder nature of PSP\textsubscript{DARP} which searches for longer AV routes (as opposed to PSP\textsubscript{CTSPAV} which only searches for mini routes). Nevertheless, regardless of the column-generation phase not converging for most instances, \emph{the DARP procedure was still able to consistently produce stronger lower bounds for the primary objective}.

The next column shows the fraction of problem instances in which the MIP for the primary objective exceeded the time limit, while the last shows the fraction of instances in which the MIP for both objectives exceeded the time limit. Both columns show that the MIP can be solved within the time limit for a majority of the problem instances, which is in stark contrast to the MIP of the CTSPAV procedure, whereby its longer time limit is exceeded for all problem instances. This difference can be attributed to the relatively easier MP\textsubscript{DARP} which just solves a set-covering problem (instead of the MP\textsubscript{CTSPAV} which solves a route-scheduling problem) and the stronger primal lower bound provided by its linear relaxation.

\begin{table}[!t]
	\caption{Optimality Gaps and Computation Times of the DARP Procedure with the Lexicographic Objective}
	\label{tab:darp_gap_time_lex}
	\resizebox{\textwidth}{!}{
		\begin{tabular}{ccccccccc}
			\noalign{\rule{0pt}{10pt}}\hline\noalign{\smallskip}
			\begin{tabular}[c]{@{}c@{}}Cluster\\location\end{tabular} &
			\begin{tabular}[c]{@{}c@{}}Depot\\config.\end{tabular} &
			\multicolumn{1}{c}{\begin{tabular}[c]{@{}c@{}}Average\\column\\count\end{tabular}} &
			\multicolumn{1}{c}{\begin{tabular}[c]{@{}c@{}}Average\\primary\\absolute\\gap\end{tabular}} &
			\multicolumn{1}{c}{\begin{tabular}[c]{@{}c@{}}Average\\primary\\optimality\\gap (\%)\end{tabular}} &
			\multicolumn{1}{c}{\begin{tabular}[c]{@{}c@{}}\% $\geq$\\primary\\colgen\\time limit\end{tabular}} &
			\multicolumn{1}{c}{\begin{tabular}[c]{@{}c@{}}\% $\geq$\\colgen\\time\\limit\end{tabular}} &
			\multicolumn{1}{c}{\begin{tabular}[c]{@{}c@{}}\% $\geq$\\primary\\MIP time\\limit\end{tabular}} &
			\multicolumn{1}{c}{\begin{tabular}[c]{@{}c@{}}\% $\geq$\\MIP\\time\\limit\end{tabular}} \\
			\noalign{\smallskip}\hline\noalign{\smallskip}
			\multirow{2}{*}{\begin{tabular}[c]{@{}c@{}}Inside\\ city\end{tabular}}  & Central & 34107 $\pm$ 5602 & 2.68 $\pm$ 0.14 & 47.9 $\pm$ 2.1 & 68.2 & 100.0 & 18.2 & 36.4 \\
			\noalign{\smallskip}\cline{2-9}\noalign{\smallskip}
			& Local   & 33911 $\pm$ 5554 & 2.55 $\pm$ 0.16 & 49.0 $\pm$ 2.7 & 68.2 & 100.0 & 9.1 & 27.3 \\
			\noalign{\smallskip}\hline\noalign{\smallskip}
			\multirow{2}{*}{\begin{tabular}[c]{@{}c@{}}Outside\\ city\end{tabular}} & Central & 6801 $\pm$ 907 & 4.81 $\pm$ 0.23 & 48.2 $\pm$ 1.5 & 95.6 & 98.5 & 14.7 & 25.0 \\
			\noalign{\smallskip}\cline{2-9}\noalign{\smallskip}
			& Local   & 6698 $\pm$ 870 & 4.79 $\pm$ 0.25 & 47.6 $\pm$ 1.6 & 95.6 & 98.5 & 14.7 & 19.1 \\
			\noalign{\smallskip}\hline\noalign{\smallskip}
		\end{tabular}
	}
\end{table}

Table \ref{tab:darp_gap_time_dist} summarizes the average optimality gaps and computations times for the DARP procedure with the distance-minimization objective. Similar to Table \ref{tab:darp_gap_time_lex}, the first three columns specify cluster locations, depot configurations, and the average number of columns generated. The following two show average optimality gaps of the final MIP results. The first uses the dual lower bound, $\ceil{z_\text{LB}}$, in the gap calculation, whereas the second uses the primal lower bound produced by the CTSPAV procedure, $z^*$. For the distance-minimization objective, the CTSPAV procedure consistently produces stronger lower bounds, therefore optimality gaps calculated using its lower bound are also smaller. The dual lower bound $\ceil{z_\text{LB}}$ is weaker in this case because the column-generation phase could not be completed within the time limit for almost all problem instances. This is further highlighted by the next column which shows the fraction of problems in which the column-generation phase could not be completed. Out of the 180 problem instances in total, only four managed to achieve convergence. For these four instances the strength of $\ceil{z_\text{LB}}$ is comparable to $z^*$ from the CTSPAV procedure. Finally, the last column shows the fraction of problem instances in which the MIP exceeded its time limit. A larger fraction of instances inside the city exceeded the time limit due to their relatively larger column counts. 

\begin{table}[!t]
	\caption{Optimality Gaps and Computation Times of the DARP Procedure with the Distance Minimization Objective}
	\label{tab:darp_gap_time_dist}
	\resizebox{\textwidth}{!}{
		\begin{tabular}{ccccccc}
			\noalign{\rule{0pt}{10pt}}\hline\noalign{\smallskip}
			\begin{tabular}[c]{@{}c@{}}Cluster\\location\end{tabular} &
			\begin{tabular}[c]{@{}c@{}}Depot\\ config.\end{tabular} &
			\multicolumn{1}{c}{\begin{tabular}[c]{@{}c@{}}Average\\column\\count\end{tabular}} &
			\multicolumn{1}{c}{\begin{tabular}[c]{@{}c@{}}Average\\optimality\\gap (\%)\end{tabular}} &
			\multicolumn{1}{c}{\begin{tabular}[c]{@{}c@{}}Average\\optimality\\gap (\%)\\(CTSPAV LB)\end{tabular}} &
			\multicolumn{1}{c}{\begin{tabular}[c]{@{}c@{}}\% $\geq$ colgen \\ time limit\end{tabular}} &
			\multicolumn{1}{c}{\begin{tabular}[c]{@{}c@{}}\% $\geq$ MIP \\ time limit\end{tabular}} \\
			\noalign{\smallskip}\hline\noalign{\smallskip}
			\multirow{2}{*}{\begin{tabular}[c]{@{}c@{}}Inside\\ city\end{tabular}}  & Central & 34286 $\pm$ 5323 & 98.0 $\pm$ 2.0 & 28.5 $\pm$ 1.7 & 100.0 & 77.3 \\
			\noalign{\smallskip}\cline{2-7}\noalign{\smallskip}
			& Local   & 30712 $\pm$ 4640 & 95.2 $\pm$ 3.8 & 36.4 $\pm$ 3.1 & 100.0 & 77.3 \\
			\noalign{\smallskip}\hline\noalign{\smallskip}
			\multirow{2}{*}{\begin{tabular}[c]{@{}c@{}}Outside\\ city\end{tabular}} & Central & 6154 $\pm$ 902 & 94.8 $\pm$ 2.6 & 28.4 $\pm$ 1.6 & 96.7 & 36.7 \\
			\noalign{\smallskip}\cline{2-7}\noalign{\smallskip}
			& Local   & 5816 $\pm$ 838 & 95.5 $\pm$ 2.6 & 44.5 $\pm$ 1.8 & 96.9 & 37.5 \\
			\noalign{\smallskip}\hline\noalign{\smallskip}
		\end{tabular}
	}
\end{table}

\end{document}